\documentclass{jams-l}
\usepackage{amsmath,amsfonts,amsthm,amssymb,amscd}
\def\classification#1{\def\@class{#1}}
\classification{\null}

\DeclareFontFamily{OT1}{rsfs}{}
\DeclareFontShape{OT1}{rsfs}{n}{it}{<-> rsfs10}{}
\DeclareMathAlphabet{\mathscr}{OT1}{rsfs}{n}{it}

\DeclareMathOperator{\SL}{SL}

\DeclareMathOperator{\rank}{rank}

\DeclareMathOperator{\Cl}{Cl}

\newcommand{\Ann}{\mathscr{A}}
\newcommand{\Q}{\mathbb{Q}}
\newcommand{\pr}{\mathfrak{p}}
\newcommand{\Z}{\mathbb{Z}}
\newcommand{\C}{\mathbb{C}}
\newcommand{\R}{\mathbb{R}}

\newcommand{\ideal}{\mathscr{I}}
\newcommand{\Norm}{N}
\newcommand{\Db}{D_{\mathrm{big}}}
\newcommand{\Ds}{D_{\mathrm{small}}}
\newcommand{\order}{\mathscr{O}}

\newtheorem{prop}{Proposition}[section]
\newtheorem{thm}[prop]{Theorem}

\newtheorem{cor}[prop]{Corollary}
\newtheorem{lem}[prop]{Lemma}
\newtheorem{defn}{Definition}
\newenvironment{Rem}{{\bf Remark.}}{}
\numberwithin{equation}{section}
\begin{document}
\title{Integral points on elliptic curves and $3$-torsion in class groups}
%\AuthorHead{H. A. Helfgott and A. Venkatesh}
\subjclass[2000]{Primary 11G05, 11R29; Secondary 14G05, 11R11.}
\keywords{Class groups, elliptic curves, integral points.}
\author{H. A. Helfgott}
\address{H. A. Helfgott, Mathematics Department, Yale University, New Haven,
CT 06520, USA}
\curraddr{H. A. Helfgott, D\'epartement de math\'ematiques et de statistique,
Universit\'e de Montr\'eal, CP 6128 succ Centre-Ville, Montr\'eal QC\; H3C 3J7, Canada.}
\author{A. Venkatesh}
\address{A. Venkatesh, Massachusetts Institute
of Technology, Department of\break \mbox{Mathematics},~Cambridge, MA 02139--4307, USA}
\curraddr{A. Venkatesh,
Courant Institute of Mathematical Sciences, New York University, NY 10012,
USA.}
\thanks{The second author was supported in part by NSF Grant DMS-0245606.}
\begin{abstract}
We give new bounds for the number of integral points on elliptic curves. 
The method may be said to interpolate between 
approaches via diophantine 
techniques (\cite{BP}, \cite{HBR}) and methods based on quasiorthogonality
in the Mordell-Weil lattice (\cite{Si4}, \cite{GS}, \cite{He}). We apply
our results to break previous bounds on the number of elliptic curves
of given conductor and the size of the $3$-torsion part of the class group
of a quadratic field. The same ideas can be used to count rational points
on curves of higher genus.
\end{abstract}
\maketitle
\section{Introduction}
We give new bounds for the number of integral points on elliptic curves.
The method may be said to interpolate between
approaches via diophantine
techniques (\cite{BP}, \cite{HBR}) and methods based on quasiorthogonality
in the Mordell-Weil lattice (\cite{Si4}, \cite{GS}, \cite{He}). We apply
our results to break previous bounds on the number of elliptic curves
of given conductor and the size of the $3$-torsion part of the class group
of a quadratic field. The same ideas can be used to count rational points
on curves of higher genus.

\subsection{Conductors and class groups}
Let $N$ be a positive integer.
We show that there are at most $O(N^{0.22377\dotsc})$ elliptic curves over
$\mathbb{Q}$ of conductor $N$.
We also prove that, for every non-zero integer $D$, at
most $O(|D|^{0.44178\dotsc})$ elements
of the class group of $\mathbb{Q}(\sqrt{D})$ are $3$-torsion.
The latter result provides the first improvement
on the trivial bound of $O(D^{1/2 + \epsilon})$, whereas
the former improves on $O(N^{1/2 + \epsilon})$, which follows 
from the said trivial bound (\cite{BS}, Thm 1). The new bound on $3$-torsion
implies that there are at most $O(|D|^{0.44178\dotsc})$
cubic extensions of $\mathbb{Q}$ with discriminant $D$ (vd. \cite{hasse},
Satz 7). These results are derived from a new method
of obtaining bounds for the number of integral
(or rational) points on curves of non-zero genus.

These questions have attracted considerable interest; see, e.g.,
\cite{duke}. A number of authors have given improved bounds either
conditionally (\cite{wong}) or in the average
(\cite{dk}, \cite{murty}, \cite{sound}).
The problems are intimately linked: the size of $3$-torsion can be 
bounded above
by the number of integral points of moderate height
on the variety $y^2 + Dz^2 = x^3$, whereas elliptic curves of given
conductor correspond to $S$-integral points on a finite collection of elliptic
curves of the form $y^2 = x^3 + C$.
The question of the size of $3$-torsion is of further interest
in view of its connection to the enumeration of cubic fields
and to upper bounds for the ranks of elliptic curves. 

The techniques in this paper are valid over an arbitrary 
number field. For example, one may show that
the number of cubic extensions of a fixed number field $K$
with prescribed discriminant $\ideal$ is $\ll \Norm \ideal^{1/2 - \rho_K}$
for some $\rho_K > 0$. This can be deduced without difficulty from
the methods of this paper: the key result about point counting, Theorem \ref{thm:astarte}, is stated over a number field.

\subsection{Points on curves}
Let $E$ be an elliptic curve defined over a number field $K$. Let $S$
be a finite set of places of $K$. We wish to bound the cardinality of the set
$E(K,S)$ of $S$-integral points on $E$.

Embed the Mordell--Weil lattice
$E(K)$, modulo torsion, into $\mathbb{R}^{\rank(E(K))}$, so that the
canonical height on $E$ is taken to the square of the Euclidean norm. Regard
$E(K,S)$ as a subset of $E(K)$. One way to bound the cardinality of
$E(K,S)$ is to exploit the fact that, in a certain sense, the points
of $E(K,S)$ tend to be apart from each other. This idea is
already present in \cite{Si4}, \cite{GS}; let us consider it in the manner
of \cite{He}, \S 4.
After some modest slicing of $E(K,S)$, we see that
any two points on the same slice are separated by almost $60^{\circ}$. We
can then apply the best available results on sphere-packing \cite{KL}
to obtain a bound on the number of elements of $E(K,S)$. 
This bound (Cor.\ \ref{cor:agobio}) improves on \cite{GS} and seems
to be the best to date. Corollary \ref{cor:newcoeur}
improves on a bound of W. Schmidt \cite{Schm}. 

A major weakness of the relatively naive method discussed thus far is that it is very sensitive
to the rank of the Mordell-Weil lattice. We have used bounds for sphere-packing problems,
and such bounds typically depend exponentially on the dimension of the ambient space. 
 This makes it difficult to apply to many natural problems, including that
of $3$-torsion in quadratic class groups, where one has a relatively poor bound on the rank
of the Mordell-Weil lattice.  In general, this problem will be particularly severe
when one is bounding the number of points on $E(K,S)$ below a certain height $h_0$,
where $h_0$ is comparable to the ``height of $E$,'' i.e., the logarithm of the largest coefficient
in a Weierstrass equation of $E$. 

Our key idea to overcome this obstacle is to exploit a certain
feature of the geometry of high-dimensional Euclidean spaces, namely, the
fact that the solutions to certain special types of packing problems 
depend relatively weakly on the dimension of the ambient space. 
More precisely, consider the question: how many vectors can one 
pack into the unit sphere on $\mathbb{R}^n$ such that the angle between any two is
$\geq \theta$? It is not difficult to see (see remark after Prop.\ \ref{prop:kl}) that one can give an upper bound {\em independent of $n$}
when $\theta > \pi/2$. We will exploit a related but considerably deeper feature, namely, 
that this phenomenon
 persists (in a much weakened form) when $\theta < \pi/2$:  for
$\theta = \pi/2 - \alpha$ the work of Kabatiansky and Levenshtein gives an upper bound, for small $\alpha$, of the form
$\exp(\alpha^2 \log(\alpha^{-1}) n)$. The critical feature here is that the constant $\alpha^2 \log(\alpha^{-1})$ depends sublinearly on $\alpha$.

We shall exploit this feature by introducing a costly type of slicing of $E(K,S)$, which allows
us to increase the angle of $60^{\circ}$, and thus
lowers the bound per slice sharply; we can see the amount of slicing
as a parameter to be optimized. This slicing is carried out as follows: we choose an auxiliary prime $p$,
and partition $E(K,S)$ into the fibers of the
 reduction map $E(K,S) \rightarrow E(\mathbb{F}_p)$; the size of $p$
is our free parameter.\footnote{The fact that this type of partitioning increases the angle is an instance of a very general phenomenon: rational points on an algebraic variety repel each other more strongly if they are forced to be $p$-adically close. This is already visible for integers: if $x,y \in \mathbb{Z}$ are distinct, one has $|x-y| \geq 1$, but if $x,y$ are congruent mod $p$ one has
$|x-y| \geq p$. }

The result obtained from $90^{\circ}+\epsilon$ is the same as what arises 
from \cite{BP}, modulo the difference between the canonical and the naive
height. (This is no coincidence; as we will see, the similarity between
the two underlying procedures runs deep.) We then show that the results
depend continuously on the angle, and that $90^{\circ}$ is a locally
suboptimal choice in the interval $\lbrack 60^{\circ},90^{\circ}\rbrack$.
Thus we will be able to make a better choice within the interval, thus
obtaining a result better than the canonical-height analogue of \cite{BP},
and, in general, better than the pure bounds as well.
It is only thus that we are able to break the $h_3(D)\ll D^{1/2}$ barrier.

The same ideas can be applied to bounding the number of 
rational points (or integer points)
up to a certain height on curves of higher genus.
This matter is discussed further in \cite{EV}, where 
it is shown how to improve in certain
contexts on the exponent $2/d$ 
occurring in the work of Heath-Brown \cite{HBR} and Elkies \cite{El}.  
We have therefore provided in the present paper only a sketch of how to extend these
methods to that case -- see Section \ref{sec:beyond}.

\subsection{Relation to other work}
The techniques known up to now for bounding integral points on elliptic curves
did not suffice to improve on the estimates $O(N^{1/2 + \epsilon})$ and
$O(D^{1/2 + \epsilon})$. Our method, like many results
in Diophantine approximation, uses the fact that
integer points that are $v$-adically close tend to repel each other. 
One may see the same underlying idea in the works of
Bombieri--Pila (\cite{BP}) and Heath-Brown \cite{HBR}; for a discussion
of the parallels between their methods and those in the present paper, see
the remark at the end of section \S \ref{subs:bipbip}.

Independently and simultaneously, L. B. Pierce has proved
a bound on $h_3(D)$ that breaks $D^{1/2}$. 
Pierce's bound is
$h_3(D) \ll D^{27/56+\epsilon}$, in general; 
for $D$ with certain divisibility properties the bound
improves to $h_3(D) \ll D^{5/12+\epsilon}$. 
The methods in
\cite{Pi} are quite different from those in the present paper; they are
based on the square sieve. 

\subsection{Acknowledgments}
We would like to thank M. Bhargava, 
A. Brumer, S. David, F. Gerth, D. Goldfeld, R. Heath-Brown, H. Iwaniec,
A. J. de Jong, L. B. Pierce, J. H. Silverman and
K. Soundararajan for their advice and encouragement. 

\section{Notation and preliminaries}
\subsection{Number fields and their places}
Let $K$ be a number field. We write $\mathscr{O}_K$ for the ring of integers
of $K$, $I_K$ for the semigroup of ideals of $\mathscr{O}_K$, $\Cl(\mathscr{O}_K)$ for
the class group of $K$,
$M_K$ for the set of all places of $K$, and $M_{K,\infty}$ for
the set of all infinite places of $K$. We write $N_{K/\mathbb{Q}} \mathfrak{a}$
for the norm of an ideal $\mathfrak{a}\in I_K$.
By a {\em prime} we will mean either a
finite place of $K$, or the prime ideal corresponding thereto. 
%Thus a sum
%or product over $p\in S$, where $S$ is 
%a set of places, is a sum or product over $v\in S$, $v$ finite.

Let $v$ be a non-archimedean place of $K$, $K_v$ the
completion of $K$ at $v$, and $p$ the prime
of $\mathbb{Q}$ below $v$. 
We denote by $v(x): K_v^{*} \rightarrow \mathbb{Z}$
the valuation, normalized as usual to be surjective, 
and we shall normalize the absolute value $|\cdot|_{v}: K_v^{*} 
\rightarrow \mathbb{R}$ 
so that it extends the usual value $|\cdot|_{p}$ of $\mathbb{Q}$.
Thus,
for $x \in K_v^{*}$, $|x|_{v} = p^{-v(x)/e_v}$, where 
$e_v$ is the ramification degree of $K_v$ over $\Q_{p}$.

Given a set of places $S\subset M_K$, we write $\mathscr{O}_{K,S}$ for the
ring of {\em $S$-integers}. An $S$-integer is an $x\in K$ with $v(x)\geq 0$
for $v\notin S\cup M_{K,\infty}$.
 We write $M(S)$ for the product of all finite places
in $S$, seen as ideals.

By $G\lbrack l\rbrack$ we mean the $l$-torsion subgroup of a group $G$.
Define $h(K) = \# \Cl(\mathscr{O}_K)$, 
$h_l(K) = \#(\Cl(\mathscr{O}_K)/ \Cl(\mathscr{O}_K)^l) =
\#(\Cl(\mathscr{O}_K)\lbrack l\rbrack)$. (By $\# A$ we
mean the cardinality of a set $A$.) The number of prime ideals dividing an
ideal $\mathfrak{a}\in I_K$ is denoted by $\omega_K(\mathfrak{a})$.
If $a\in \mathbb{Z}$, we may write $\omega(a)$ instead of 
$\omega_K((a))$.

Given a place $v$ of $K$, we let
\begin{equation}\label{eq:duchess}
\gamma_v = \begin{cases} 1 &\text{if $v$ is infinite}\\
0 &\text{if $v$ is finite.} \end{cases}\end{equation}
We further set $d_v = [K_v:\Q_p]$ where $p$ is the place
of $\Q$ below $K$; in particular, $d_v = 2$ or $1$
when $v$ is complex or real, respectively.

If $R$ is an integral domain with quotient field $K$,
and $M$ is an $R$-module, we write
$\rank_R(M)$ for the dimension of $M\otimes_R K$ over $K$. 
 
For
every $r\in \mathbb{R}$, we define
\[\log^+ r  = \log(\max(r,1)) .\]
Given $x\in K$, we define its {\em height}
\begin{equation}\label{eqn:htrel} h_K(x) = \sum_{v\in M_K} d_v
\log^+(|x|_v)\end{equation}
and its {\em absolute height}
\begin{equation} \label{eqn:htabs}
h(x) = \frac{1}{\lbrack K : \mathbb{Q}\rbrack} h_K(x). \end{equation}

\subsection{Elliptic curves}
Let $E$ be an elliptic curve over a number field
$K$. Given a field $L\supset K$,
we use $E(L)$ to denote the set of $L$-valued points of $E$. (Thus
$E(\mathbb{Q})$ is the set of rational points of an elliptic curve
defined over $\mathbb{Q}$.) 
We take $E$ to be given by a Weierstrass equation
\begin{equation} \label{eq:weierstrass}
E: y^2 + a_1 x y + a_3 y = x^3 + a_2 x^2 + a_4 x + a_6,\end{equation}
where $a_1,\dotsc,a_6\in \mathscr{O}_K$.
By $j(E)$ (resp. $\Delta(E)$) we mean the $j$-invariant
(resp. discriminant) of (\ref{eq:weierstrass}). 
We write $x(P)$, $y(P)$ for the $x$- and $y$-coordinates of a point
$P\in E(\overline{K})$ other than the origin.
Given a set of places $S\subset M_L$,
we denote by $E(L,S)$ the set of {\em $S$-integral points}, i.e.,
points with $S$-integral coordinates:
\[E(L,S) = \{P\in E(L)\backslash \{0\} : x(P), y(P) \in \mathscr{O}_{L,S}\} .\]

As is usual, we write $\hat{h}$ for the {\em canonical height} on $E$, defined
on all points of $E(\overline{K})$. 
The canonical height $\hat{h}$ is a positive definite quadratic form\footnote{where ``positive definite'' is taken to mean 
`mapping {\em non-torsion} elements to
positive numbers.''}
on the abelian group $E(\overline{K})$, or, by restriction, on
$E(K)$. It lets itself be expressed as a sum of {\em local height functions}
$\lambda_v:E(K_v)\to \mathbb{R}$, as follows:
\[\hat{h}(P) = \frac{1}{\lbrack K : \mathbb{Q}\rbrack} \hat{h}_K(P) = 
 \frac{1}{\lbrack K: \mathbb{Q}\rbrack} \sum_{v\in M_K} d_v \lambda_v(P) .\]

Local height functions are canonically defined up to an additive constant;
we follow the conventions in \cite{La}, Ch.\ VI, and \cite{Si2}, Ch.\ VI,
which make local heights independent of the model. Note that 
$\lambda_v(P) = \lambda_v(-P)$ for every place $v$, and, in particular,
$\lambda_v(P_1 - P_2) = \lambda_v(P_2-P_1)$ for any place $v$ and any
$P_1, P_2\in E(K_v)$.

We recall that an elliptic curve $E$ over a nonarchimedean local field $K$ is said to have {\em potentially 
good reduction} if it admits a model with good reduction
in some extension of $K$. We say that $E$ has
{\em potentially
multiplicative reduction}
 if it does not have potentially good reduction; this occurs precisely when the $j$-invariant of $E$ is not integral. See \cite[Chapter VII]{Si}.  

\section{Integral points on elliptic curves}
\subsection{Uniform quasi-orthogonality}
Integral points on elliptic curves tend to repel each other; so do
rational points on curves of higher genus. A classical formulation of
the latter fact is due to Mumford \cite{Mu}; the former phenomenon can be 
seen to surface in \cite{Si4} and \cite{GS}. In order to go
further, however, we must quantify this repulsion in a fashion that is more uniform and
more flexible than those available up to date.

As in \cite{GS}, we will use local heights.
Roughly speaking, we wish to establish a result
of the form $\lambda_v(P-Q) \geq \mathrm{min}(\lambda_v(P),
\lambda_v(Q))$.
Although this is not quite true at places of bad reduction or at the
archimedean
places, it is true if we subdivide $E(K_v)$ into a fairly small number of
slices and ask that $P,Q$ lie in the same slice; see
Lemmas \ref{lem:potgood}--\ref{lem:tetris}. 
One feature of these Lemmas is that they provide
somewhat sharper results in the region where $\lambda_v(P) \leq 0$ than
elsewhere; 
this will eventually be significant in dealing with points of small
(global) height.
We can then prove the
quasi-orthogonality result in Prop.\ \ref{prop:affreuse}. 
In words, it asserts: integral points
are quite far apart from each other in the Mordell-Weil lattice, and, moreover, forcing two integral points to be
congruent modulo some ideal of $\order_K$ forces them even further apart in the Mordell-Weil lattice. 

It should be remarked that if one is willing to accept an extra factor of size
about $(1+m)^m$
in Thm \ref{thm:astarte}, where $m$ is the number of places of 
potentially multiplicative reduction, the proofs that follow can be considerably simplified. 
In this context, note that
$y^2 = x^3 + D$ has in fact $m=0$, so this weaker version
would suffice for the applications in Section \ref{sec:ghatyar}.
Indeed, for the applications of Section \ref{sec:ghatyar}, it is not difficult to
avoid local heights completely: since we deal with the curves $y^2 = x^3 +
D$, one may use the fact that they are all twists of $y^2 = x^3 + 1$
to prove the required special cases of Prop.\ \ref{prop:affreuse} and Thm.\
\ref{thm:astarte} in an elementary fashion (cf.\ \cite{He}, Lem.\ 4.16).

\begin{lem}\label{lem:potgood}
 Let $E$ be an elliptic curve over a non-archimedean local field $K_v$
 with potentially good reduction.
 Let $P_1, P_2\in E(K_v)$ be two distinct points. Then
 \[\lambda_v(P_1-P_2)\geq \min(\lambda_v(P_1),\lambda_v(P_2)) .\]
 \end{lem}
 \begin{proof}
 Pass to an extension $L_w$ of $K_v$ on which $E$ acquires good
 reduction.
 Choose a Weierstrass equation for $E$ over $L_w$ such that $v(\Delta)=0$.
 Then $\lambda_v(P) = \lambda_w(P) = \frac{1}{2} \log^+(|x(P)|_w)$.
The statement follows therefrom by direct computation. (Alternatively, use the
definition of the local height in terms of the canonical filtration.)
%The statement then follows from
%$|x(P_1-P_2)|_w\leq \min(|x(P_1)|_w,|x(P_2)|_w)$.
 \end{proof}

\begin{lem}\label{lem:potmult}
 Let $E$ be an elliptic curve over a non-archimedean local field $K_v$
 with potentially multiplicative reduction. Then, for 
any sufficiently small $\epsilon>0$, there is a partition
 \begin{equation}\label{eq:marmite}
 E(K_v) = W_{v,0} \cup W_{v,1} \cup\dotsb \cup W_{v,n_v},\;\;
 n_v\ll |\log \epsilon|,\end{equation}
 such that for any two distinct points $P_1, P_2\in W_{v,0}$,
 \[\lambda(P_1-P_2)\geq
 \min(\lambda(P_1),\lambda(P_2))\;\;\:\text{and}\;\;\;
 \lambda(P_1),\lambda(P_2)\geq 0,\]
 and for any two distinct points $P_1, P_2\in W_{v,j}$, $1\leq j\leq n_v$,
 \[\begin{aligned}
\lambda(P_1 - P_2) &\geq (1 - \epsilon) \max(\lambda(P_1),\lambda(P_2)),\\
\lambda(P_1 - P_2) &\geq (1 - 2 \epsilon) \max(\lambda(P_1),\lambda(P_2)) .
\end{aligned}\]
 The implied constant is absolute.
 \end{lem}
 \begin{proof}
 %First of all, we cut and paste the last paragraph before Lemma 3.1, sans
 % the first line and a half:
 The elliptic curve $E$ is isomorphic, over an algebraic closure $\overline{K_v}$,
 to a Tate curve $E_{q}$ for some $q \in K_v^*$ satisfying
$v(q) = - v(j)$, where $j = j(E)$ is the $j$-invariant of $E$; see
\cite{Si2}, Ch.\ V.
  
 There is a natural composition
 \[\beta_v:E(K_v) \to E(L_w) \to
 E(L_w)/E_0(L_w) \stackrel{\alpha}{\rightarrow}
  \mathbb{R}/\mathbb{Z} \rightarrow \lbrack 0,1),\]
  where
  $L_w/K_v$ is the minimal extension such that $E$ acquires
  split multiplicative reduction over $L_w$, and the map $\alpha$ is
  given by $\alpha(t) = v(t)/v(q)$ on the Tate curve. 
  (See \cite{La}, pp.\ 68--69; cf.\
  \cite{GS}, p.\ 270.)
  %Notice that it is important to have [0,1) rather than [-1/2,1/2) above -
  %see La, end of p.\ 69, if you want to know the reason. We will be able
  %to work with [-1/2,1/2) when we do the archimedean places.
  %Now we proceed with the proof as it was before:
  For every $P\in E(L_w)$,
  \[\lambda(P) = - \frac{1}{2} B_2(\beta_v(P)) \log |q|_w - \iota(P) \log
  |\pi|_w,\]
  where $w$ is an extension of $v$ as above,
  $\pi$ is a uniformizer of $w$, $B_2(t) = t^2 - t + 1/6$, $\iota(P) = 0$
  when $P\notin E_0(L_w)$, and, if $P\in E_0(L_w)$, $\iota(P)$
is the largest integer $\iota$ such that $P\in E_{\iota}(L_w)$, where
\[E(L_w) \supset E_0(L_w) \supset E_1(L_w) \supset \dotsb\]
is the canonical filtration of $E(L_w)$. (See \cite{La}, pp.\ 68--70; cf.\
\cite{GS}, p.\ 270.)
 Clearly
 \[\iota(P_1-P_2) \geq \min(\iota(P_1),\iota(P_2)) .\]
 If $\beta_v(P_1) = \beta_v(P_2) = 0$, it follows that $\lambda(P_1 - P_2)
 \geq \min(\lambda(P_1), \lambda(P_2))$.
 %Now we get started with the splitting into regions W_j.
 Define, then, $W_{v,0} = \{P\in E(K_v) : \beta_v(P)=0\}$.
 It remains to partition $\{P\in E(K_v) : \beta_v(P)\ne 0\}$.
 %We cut and paste from the proof of what used to be Lemma 3.2 (now Lem.\
 %3.3)
 Partition $(0,1/2\rbrack$ into sets $U_0, U_1,\dotsc, U_m$,
 where $m = \lceil \log_{3/2} (6/\epsilon)\rceil$:
 \[\begin{aligned}U_0 &= (0, \epsilon/12\rbrack,\; U_m = ((3/2)^{m-1}
 \epsilon/12,
 1/2\rbrack,\\
 U_j &= ((3/2)^{j-1} \epsilon/12, (3/2)^j \epsilon/12\rbrack,\;\;\; 1\leq j<m
 .
 \end{aligned}\]

Note that $B_2(t)$ is decreasing on $t \in [0,1/2]$.
Suppose $t_1,t_2$ both belong to $U_j$ and $t_1 \geq t_2$.  If $j=0$,
we have 
\begin{eqnarray} \label{eq:jeq0} B_2(t_1-t_2) \geq B_2(\epsilon/12)
\geq 1/6 - \epsilon/12 =  (1 - \epsilon) B_2(0) + \epsilon/12
\\ \nonumber > (1-2\epsilon) B_2(0) + \epsilon/12
\end{eqnarray}
If $j \geq 1$ then $u \leq t_2 \leq t_1 \leq 3u/2$ 
where $u = (3/2)^{j-1} \epsilon/12$.  Then:
\begin{equation} \label{eq:jeq1}\begin{aligned}
B_2(t_1 - t_2) &\geq B_2(u/2),\\
B_2(u/2)  \geq (1-\epsilon)B_2(u) + \epsilon/12,\;\;\;\;
& B_2(u/2)  \geq (1-2 \epsilon)B_2(u) + \epsilon/12 .
\end{aligned}\end{equation}

(The last two inequalities are proved in two
cases according to whether $B_2(u) \geq 1/12$ or
$B_2(u) < 1/12$. In the former case
we have $B_2(u/2) \geq B_2(u) = (1-\epsilon)B_2(u) + \epsilon B_2(u)
\geq (1-\epsilon) B_2(u) + \epsilon/12 > (1 - 2 \epsilon) B_2(u)
+ \epsilon/12$. In the latter case, 
$u > 1/11$ and $B_2(u/2) - B_2(u)  > 1/30$; 
in particular $B_2(u/2) - (1-\epsilon) B_2(u) - \epsilon/12 \geq
1/30 + \epsilon B_2(u) - \epsilon/12$. Since $B_2(u) \geq -1/12$,
it follows that $B_2(u/2) - (1-\epsilon) B_2(u) - \epsilon/12
\geq 1/30 - \epsilon/6 \geq 0$, where we assume $\epsilon < 1/5$;
similarly, $B_2(u/2) - (1 - 2 \epsilon) B_2(u) - \epsilon/12 \geq
1/30 - \epsilon/4 \geq 0$, where we assume $\epsilon < 2/15$.)

 Combining (\ref{eq:jeq0}) and (\ref{eq:jeq1})
we obtain\footnote{The term $\epsilon/12$ in
 the displayed equation will be used in the proof of Lem.\ \ref{lem:tetris}.}
 \begin{equation} \label{eq:boundm} \begin{aligned}
 B_2(t_1-t_2) &\geq (1 - \epsilon)
 \max_{j=1,2} B_2(t_j) + \epsilon/12 \\
 B_2(t_1-t_2) &\geq (1 - 2\epsilon)
  \max_{j=1,2} B_2(t_j) + \epsilon/12 \end{aligned} \end{equation}
 for all $t_1, t_2\in U_j$, $0\leq j\leq m$, where $t_1 \geq t_2$.  Define
 \[\begin{aligned}
 W_{v,2 j + 1} &= \{P\in E(K_v) : \beta_v(P) \in U_j\},\\
 W_{v,2 j + 2} &= \{P\in E(K_v) : \beta_v(-P) \in U_j,\,
 \beta_v(P)\ne 1/2\} \end{aligned}\]
 for $0\leq j\leq m$. We set $n_v = 2 m + 2$ and are done.
 \end{proof}
 %In other words, we learned a couple of tricks from the proof of
 %the following Lemma.
                                            
\begin{lem}\label{lem:tetris}
Let $E$ be an elliptic curve over $\mathbb{C}$. 
Then, for any sufficiently small $\epsilon>0$, there is a partition
\begin{equation}\label{eq:ago}
E(\mathbb{C}) = W_0 \cup W_1 \cup \dotsb \cup W_n,\;\;
n\ll \epsilon^{-2} |\log \epsilon|^2,\end{equation}
such that for any two distinct points $P_1, P_2\in W_{j}$,
$0\leq j\leq 5$,
\begin{equation} \label{eq:required2}
\lambda(P_1-P_2)\geq (1 - \epsilon)
\min(\lambda(P_1),\lambda(P_2))\;\;\:\text{and}\;\;\;
\lambda(P_1),\lambda(P_2)\geq 0,\end{equation}
and for any two distinct points $P_1, P_2\in W_{j}$, $6\leq j\leq n$,
\begin{equation} \label{eq:required}\begin{aligned}
\lambda(P_1 - P_2) &\geq (1 - \epsilon) \max(\lambda(P_1),\lambda(P_2)),\\
\lambda(P_1 - P_2) &\geq (1 - 2 \epsilon) \max(\lambda(P_1),\lambda(P_2)) .
\end{aligned}\end{equation}
The implied constant in (\ref{eq:ago}) is absolute.
\end{lem}
\begin{proof}
There is an isomorphism $E(\mathbb{C}) \stackrel{u}{\rightarrow}
 \mathbb{C} / (\mathbb{Z} + \tau
 \mathbb{Z})$ for some $\tau$ in the usual fundamental domain of
 $\SL_2(\mathbb{Z})\backslash \mathbb{H}$. Note especially that
 $\Im(\tau) \geq \sqrt{3}/2$. Write
   $u(P) = u_{P,1} + \tau u_{P,2}$,
   $u_{P,1}, u_{P,2} \in \left(-\frac{1}{2},\frac{1}{2}\right\rbrack$.
    Define $q = e^{2\pi i \tau}$,
    $q_u(P) = e^{2\pi i u_P}$ --  note that $|q| \leq e^{-\pi \sqrt{3}}$.
    %The following paragraph is much abbreviated, as it is already present
    %in the main in the proof of the previous Lemma.
    The local height is given by
    \begin{equation}\label{eq:monst}
    \lambda(P) = -\frac{1}{2} B_2(u_{P,2}) \log |q| - \log |g_0(q_u(P))| ,
    \end{equation}
    where
   \begin{equation} \label{eq:strom1}
    g_0(t) = (t-1) \prod_{n=1}^{\infty} (1 - q^n t) (1 - q^n t^{-1})
        \end{equation}
    and $B_2(t)$ is as in the proof of Lem.\ \ref{lem:potmult}. 
   (See, e.g. \cite{Si2}, Ch.\ VI, Thm 3.4.) We partition
    $[0,1/2\rbrack$ into sets $U_0, U_1,\dotsc, U_m$ as in the same
    proof. For the present proof we adjoin $0$ to $U_0$, since Lem.\
    \ref{lem:potmult} partitions only $(0,1/2)$. 

Whenever $t_1, t_2$, with $t_2 \leq t_1$, belong to the same set $U_j$ 
one obtains
\begin{equation}\label{eq:tutankhamen}\begin{aligned}
B_2(t_1 - t_2) &\geq (1-\epsilon) \max_{j=1,2} B_2(t_j) + \epsilon/12,\\
B_2(t_1 - t_2) &\geq (1 - 2 \epsilon) \max_{j=1,2} B_2(t_j) + \epsilon/12 
\end{aligned}\end{equation}
as in (\ref{eq:boundm}), where we assume $\epsilon<2/15$.

Let $P_1, P_2\in E(\mathbb{C})$.
Since $\lambda(P) = \lambda(-P)$, we may assume
 $0\leq u_{P_2,2}\leq u_{P_1,2}\leq 1/2$ without loss of generality.
 Let $\Ann$ be the annulus $\{z: |q|^{1/2} \leq |z| \leq
 1\}$; thus whenever $0 \leq u_{P,2} \leq 1/2$ we have $q_u(P) \in \Ann$.

 In view of (\ref{eq:tutankhamen}) and
$- \frac{1}{2} \log |q| \geq \frac{\pi \sqrt{3}}{2} > 2$,
 it will suffice (see (\ref{eq:strom3})) to partition $\mathscr{A}$ into sets
  $V_0,V_1,\dotsc,V_{m'}$ such that
  \begin{equation}\label{eq:varlam}\begin{aligned}
  -\log |g_0(q_u(P_1-P_2))| &\geq (1-\epsilon) \min_{j=1,2} (- \log
  |g_0(q_u(P_j))|) - \epsilon/6, \\
0 \leq u_{P_2,2}, u_{P_1,2} &\leq \epsilon/12,\;\;\;\;\;
\lambda(P_1), \lambda(P_2) \geq 0
  \end{aligned}\end{equation}
if  $u_{P_2,2} \leq u_{P_1,2}$ and $q_u(P_1), q_u(P_2)\in V_j$  for some $j\in \{0,1,2\}$, and
\begin{equation}\label{eq:recall}\begin{aligned}
-\log |g_0(q_u(P_1-P_2))| &\geq (1-\epsilon) \max_{j=1,2} (- \log
|g_0(q_u(P_j))|) - \epsilon/6,\\
-\log |g_0(q_u(P_1-P_2))| &\geq (1-2\epsilon) \max_{j=1,2} (- \log
|g_0(q_u(P_j))|) - \epsilon/6
\end{aligned}\end{equation}
if $u_{P_2,2} \leq u_{P_1,2}$ and $q_u(P_1), q_u(P_2)\in V_j$ for some $j\in \{3,4,\dotsc,m'\}$.

%(To see that a partition $V_0, \dots, V_m'$
%which satisfies (\ref{eq:varlam}) and (\ref{eq:recall}) will suffice
%for the statement of the Lemma, we proceed as follows.
%We partition $E(\C)$ into the least common refinement of the partitions
%$q_u^{-1}(V_{*})$ and $u^{-1}([0,1] \times U_{*})$. 
%
%First suppose $P_1, 
%P_2 \in V_j$ for $j \in \{3,4, \dotsc,m'\}$
%satisfy $u_{P_2,2} \leq u_{P_1, 2}$ and $u_{P_1,2}, u_{P_2,2} \in U_r$,
%some $0 \leq r \leq m$. 
%Taking (\ref{eq:recall}) and (\ref{eq:boundm}) 
%together shows that, for $k \in \{1,2\}$, we
%have
%\begin{eqnarray*}\lambda(P_1-P_2)
%\geq (1-\epsilon) B_2(|u_{P_k,2}|) (\frac{-\log|q|}{2} )
%+ \frac{\epsilon}{12} \frac{-\log|q|}{2} \\ + 
%(1-\epsilon)(-\log|g_0(q_u(P_k))|) - \epsilon/40
% \geq (1-\epsilon) \lambda(P_k) + \epsilon(\frac{-\log|q|}{24}
%- \frac{1}{40}).\end{eqnarray*}
%This implies (\ref{eq:required}) since $-\log|q|/24 \geq 1/40$.
%
%Now suppose that $P_1, P_2 \in V_0$
%satisfy $0 \leq u_{P_2,2} \leq u_{P_1,2} \leq \epsilon/6$.
%Let $k \in \{1,2\}$ be such that $-\log|g(q_u(P_k)| 
%= \min_{j \in \{1,2\}} -\log|g(q_u(P_j))|$. 
%(\ref{eq:boundm}) and (\ref{eq:varlam}) now imply
%that 
%\begin{eqnarray*}\lambda(P_1 - P_2) \geq 
%\frac{-\log|q|}{2} (1- \epsilon) B_2(u_{P_k,2})
%+ \frac{\epsilon}{12} \frac{-\log|q|}{2} \\ + 
%(1-\epsilon) (-\log|g_0(q_u(P_k))|)  - \epsilon/12
%\\ = (1-\epsilon) \lambda(P_k) + \frac{\epsilon}{12}(-\frac{\log|q|}{2} -
%1)
%\end{eqnarray*}
%This implies (\ref{eq:required2})
%since $\frac{\epsilon}{12}(\frac{-\log|q|}{2} -1) \geq 0$. )
{\em The region near $z=1$.}
Let $\Db = \{\delta\in \mathbb{C}: |\delta|\leq 1/2\},
\Ds = \{\delta \in \mathbb{C}: |\delta| \leq 1/8\}$.
Note $(1 + \Ds) . (1+ \Ds)^{-1} \subset 1+\Db$.

For $t = 1 + \delta$, $\delta\in \Db$,
\begin{equation}\label{eq:nonsense}\begin{aligned}
- \log |g_0(t)| &= - \log |t-1| - \sum_{n=1}^{\infty}
\log |1 - q^n (t + t^{-1})
 + q^{2 n}|\\
 &= -\log |\delta| - 2 \sum_{n=1}^{\infty} \log |1 - q^n| + O(|\delta|^2) .
 \end{aligned}\end{equation}
 (Here we use the fact that $|q| \leq e^{-\pi \sqrt{3}}$.)
 Thus, for any $\delta_1, \delta_2\in \Ds$ with $\arg(\delta_1/\delta_2)\in
 \lbrack -\frac{\pi}{3},\frac{\pi}{3}\rbrack$,
 \begin{equation}\label{eq:polen}\begin{aligned}
 -\log \left|g_0\left(\frac{1+\delta_1}{1 + \delta_2}\right)\right| &=
 - \log |\delta_1 - \delta_2| -
 2 \sum_{n=1}^{\infty} \log |1 - q^n| + O(\max(|\delta_1|,|\delta_2|))\\
 &\geq \min_{j=1,2} (-\log |g_0(1 + \delta_j)|) + O(\max_{j=1,2} |\delta_j|)
 .\end{aligned}\end{equation}
 We can thus define the sets
 \[V_k = \left\{ z\in \mathscr{A} :
 |1 - z|\leq \kappa_0 \epsilon,\,
 \arg(1-z) \in \left\lbrack -\frac{\pi}{2} + \frac{\pi}{3} k,
 -\frac{\pi}{2} + \frac{\pi}{3} (k + 1) \right\rbrack \right\},\]
 where $k=0,1,2$ and $\kappa_0$ is small enough so that (a)
 $O(\max_{j=1,2} |\delta_j|)$ in (\ref{eq:polen}) is 
less than $\epsilon/6$ in absolute value when
 $|\delta_j| \leq \kappa_0 \epsilon$, 
(b) 
$|\log |1 - \kappa_0 \epsilon| |/
  (\pi \sqrt{3}) \leq \epsilon/12$, and
 (c) $- \log |g_0(1+\delta)|\geq 0$, for any $q$, whenever $|\delta|\leq
 \kappa_0 \epsilon$.
 The conditions in (\ref{eq:varlam}) are then satisfied.

{\em The region near z=0.} We will partition the region
$\{z\in \mathscr{A} :
|z|\leq \kappa \epsilon\}$ for some constant $\kappa$.

%Consider the functions
%f_{+}(x) = \prod_{n=-\infty}^{\infty} (1 + x^{2 n +1})$,
%$f_{-}(x) = \prod_{n=-\infty}^{\infty} (1 - x^{2 n + 1})$. There is
%an absolute constant $\kappa_0$ such that $f_+(x)\leq e^{\epsilon/36}$ and
%$f_-(x) > e^{-\epsilon/36}$ for $0\leq x\leq \kappa_0 \epsilon$. Thus

For $t \in \mathscr{A}$ we have the bounds
\begin{equation} \label{eq:strom2}
\begin{aligned}
\prod_{n=1}^{\infty} (1-|q|^{n-1/2} - |q|^{n+1/2} - |q|^{2 n}) &\leq 
\left|\prod_{n=1}^{\infty} (1-q^n t) (1-q^n t^{-1})\right|,\\
\left|\prod_{n=1}^{\infty} (1-q^n t) (1-q^n t^{-1})\right| &\leq
\prod_{n=1}^{\infty} (1+|q|^{n+1/2}) (1+|q|^{n-1/2})
.\end{aligned}\end{equation}
In particular, there is an absolute constant $\kappa_1$
such that, if $t \in \mathscr{A}$ and $|q|^{1/2} \leq \kappa_1 \epsilon$, 
\begin{equation} \label{eqn:froggy}
e^{-\epsilon/18} \leq \left| \prod_{n=1}^{\infty} (1 - q^n t) (1 - q^n
t^{-1})\right| \leq e^{\epsilon/18} .\end{equation}

 We will eventually choose $\kappa \leq \kappa_1$, 
so that if
$|q|^{1/2} > \kappa_1 \epsilon$, then $|z|> \kappa_1 \epsilon$ for all
$z\in \mathscr{A}$ and
the set $\{t \in \mathscr{A}: |t| < \kappa \epsilon\}$ is empty. 
We may therefore assume that 
$|q|^{1/2} \leq \kappa_1 \epsilon$
and that (\ref{eqn:froggy}) holds. Now, for any $t\in \mathscr{A}$ such that
$e^{-\epsilon/18} \leq |t-1| \leq e^{\epsilon/18}$,
\begin{equation} \label{eq:cojones} |- \log |g_0(t)|| =
\left|- \log |t-1| - \log \left|
 \prod_{n=1}^{\infty} (1 - q^n t) (1 - q^n t^{-1})\right|\right| \leq 
 \epsilon/9
 .\end{equation}

 For $k=1,2,\dotsc,6$, let
 \[V_{k+2} = \left\{z\in \mathscr{A} : |z|\leq \kappa_2 \epsilon,\,
 \arg(z) \in \left\lbrack \frac{(k-1) \pi}{3}, \frac{k
 \pi}{3}\right)\right\},\]
 where $\kappa_2$ is an absolute constant such that
 $e^{-\epsilon/18} \leq |z-1|\leq e^{\epsilon/18}$ for $|z|\leq \kappa_2 \epsilon$. 
 Suppose $P_1, P_2 \in E(\mathbb{C})$ are such that $0 \leq
u_{P_2,2} \leq u_{P_1,2} \leq 1/2$
and $q_u(P_1), q_u(P_2) \in V_{k+2}$ for some $1 \leq k \leq 6$. Then
$q_u(P_1 - P_2) $ belongs to $\mathscr{A}$ and satisfies
$|q_u(P_1 - P_2) - 1| \leq 1$. Then (\ref{eqn:froggy})
shows that $-\log|g(q_u(P_1-P_2))| \geq -\epsilon/18$. 
Combining this with (\ref{eq:cojones}), we obtain:
 \begin{equation}\begin{aligned}
-\log |g_0(q_u(P_1-P_2))| &\geq 
 (1 - \epsilon) \max_{j=1,2} (- \log |g_0(q_u(P_j))|)  - \epsilon/6,\\
-\log |g_0(q_u(P_1-P_2))| &\geq 
 (1 - 2 \epsilon) \max_{j=1,2} (- \log |g_0(q_u(P_j))|)  - \epsilon/6 
\end{aligned}\end{equation}
  for $P_1, P_2\in E(\mathbb{C})$ with $u(P_1), u(P_2)\in V_{k+2}$. We set
  $\kappa = \min(\kappa_1, \kappa_2)$ and are done.

{\em The remaining region.} It remains to partition the region
$R = \{z\in \mathscr{A} : |z|> \kappa \epsilon, |1-z|> \kappa_0 \epsilon\}$.

By virtue of 
(\ref{eq:strom1}) and (\ref{eq:strom2}), if $q_u(P) \in \mathscr{A}$,
then $-\log|g(q_u(P))|$ differs from
$-\log|q_u(P)-1|$ by an absolutely bounded constant. 
In particular, there are absolute constants $c_1, c_2,
c_3$ such 
that, for any $c < 1$,
\begin{equation} \label{eq:bound1}
-\log|g_0(q_u(P'))| \geq -\log(c \epsilon) + c_1\end{equation}
whenever $q_u(P') \in \mathscr{A}, |q_u(P')-1| \leq c \epsilon$, and
\begin{equation} \label{eq:bound2}
c_3 \leq -\log|g_0(q_u(P))| \leq -\log(\epsilon) + c_2\end{equation}
for all $P \in R$. By
(\ref{eq:bound1}) and (\ref{eq:bound2}),
we may choose a sufficiently small (absolute) constant
$c$ such that
\begin{equation} \label{eq:donaldtrump}\begin{aligned}
- \log|g_0(q_u(P'))|
&\geq -(1- \epsilon) \log|g_0(q_u(P))|,\\
- \log|g_0(q_u(P'))|
&\geq -(1-2 \epsilon) \log|g_0(q_u(P))|
\end{aligned}
\end{equation} for all $P$, $P'$ with $q_u(P)\in R$,
$q_u(P')\in \mathscr{A}$, $|q_u(P') - 1| < c \epsilon$.

Now, for any $P_1$,
$P_2$ with $q_u(P_1), q_u(P_2)\in R$,
$0\leq u_{P_2,2}\leq u_{P_1,2}\leq 1/2$ and
\[|\Re \log \frac{q_u(P_1)}{q_u(P_2)}|,
|\Im \log \frac{q_u(P_1)}{q_u(P_2)}| \leq \frac{c \epsilon}{2},\]
we have  $|q_u(P_1 - P_2)-1|< c \epsilon$. Hence it is enough to partition
$\log(R)$ into squares of side $c \epsilon/2$.
Since $\log(R)$ is contained in the
rectangle $\lbrack \log(\kappa \epsilon), 0\rbrack \times \lbrack -\pi,\pi
\rbrack$, there are $O(\epsilon^{-2} |\log \epsilon|)$ such squares. Their
images under exp partition $R$ into sets $V_9,V_{10},\dotsc,V_{m'}$,
with $m'\ll \epsilon^{-2} |\log \epsilon|$.
We have
$q_u(P_1-P_2) \in \mathscr{A},
|q_u(P_1-P_2) - 1| \leq c \epsilon$
if $P_1, P_2 \in V_k, 0 \leq u_{P_2,2} \leq u_{P_1,2} \leq 1/2$ 
for some $9 \leq k \leq m'$.
By (\ref{eq:donaldtrump}), we may conclude that, 
for $P_1, P_2\in E(\mathbb{C})$ with $u(P_1), u(P_2)\in V_k$, $9\leq k\leq
m'$, we have:
\[\begin{aligned}
-\log|g_0(q_u(P_1-P_2))| &\geq (1- \epsilon) \max_{j=1,2}(-
\log|g_0(q_u(P_j))|),\\
-\log|g_0(q_u(P_1-P_2))| &\geq (1-2 \epsilon) \max_{j=1,2}(-
\log|g_0(q_u(P_j))|),
\end{aligned}\]
which certainly imply (\ref{eq:recall}). 

{\em Conclusion.}
Let $u_2: E(\mathbb{C}) \rightarrow (-1/2,1/2]$ be the map
$P \mapsto u_{P,2}$. 
We partition $E(\mathbb{C})$ into the sets
\begin{equation} \label{eq:strom3}
u_2^{-1}(U_i) \cap q_u^{-1}(V_k),\;\;\;
u_2^{-1}(\{x\in - U_i: x\ne 0, -1/2\}) \cap (-q_u^{-1}(V_k)),\\
\end{equation}
where $0\leq i\leq m$, $0\leq k\leq m'$.
The inequality
$0 \leq u_{P_2,2}, u_{P_1,2} \leq \epsilon/12$ in (\ref{eq:varlam}) 
ensures that $u_2^{-1}( U_i) \cap q_u^{-1}(V_k) = \emptyset$
for $k=0,1,2$, $i\ne 0$. We define $W_0, \dots, W_5$
to be the sets in (\ref{eq:strom3}) arising from $0 \leq k \leq 2, \, i =0$,
and $W_6, W_7, \dots, W_m$ to be the other non-empty sets in
(\ref{eq:strom3}).

\end{proof}

%When invoking Lemmas \ref{lem:potmult} and \ref{lem:tetris}
%in what follows, we shall use the notation $W_i(\epsilon)$
%when it is necessary to make clear the dependence of
%the sets $W_i$ on $\epsilon$. 

\begin{prop}\label{prop:affreuse}
Let $E$ be an elliptic curve over a number field $K$,
given by a Weierstrass equation (\ref{eq:weierstrass}). 
Let $S$ be a finite set of places of $K$,
including all infinite places and all primes dividing
the discriminant of $E$. 

Let $P_1, P_2 \in E(K,S)$ be two distinct $S$-integral points.
Let
$\epsilon$ be sufficiently small. Assume that $P_1$ and $P_2$ belong to the
same set $W_{v,i}$ for every infinite place $v$ and 
every place $v$ of
potentially multiplicative reduction (see (\ref{eq:marmite}),
(\ref{eq:ago})).
Furthermore, suppose that
\begin{equation}\label{eq:frenchfries}
\sum_{v \in  T} d_v |\lambda_v(P_1) - \lambda_v(P_2)| \leq
\epsilon \max_{j=1,2} \sum_{v\in  T} d_v \lambda_v(P_j) ,\end{equation}
where $T=\{v\in S : \lambda_v(P_1) \geq 0 \mbox{ and } \lambda_v(P_2)\geq 0\}$.
Let $\ideal$ be an ideal of $\mathscr{O}_K$ not divisible by any
primes in $S$.
Assume that $P_1$ and $P_2$ have the same reduction\footnote{
Under the stated assumptions on $E$ and $\ideal$,
there is a well-defined reduction map $E(K) \rightarrow
E(\mathscr{O}_K/\ideal)$. 
}
modulo $\ideal$.
Then
\[\hat{h}(P_1-P_2)\geq (1-2 \epsilon) \max_{j=1,2} \hat{h}(P_j)
+ \frac{\log(\Norm \ideal)}{\lbrack K : \mathbb{Q}\rbrack} . \]
\end{prop}
\begin{proof}
For every finite place $v$ of good reduction, $\lambda_v(P)\geq 0$
(by e.g. \cite{La}, Thm.\ VI.4.3, or \cite{Si2}, Thm.\ VI.4.1).
Hence \[\hat{h}_K(P_1 - P_2) \geq \sum_{v\in S} d_v \lambda_v(P_1 - P_2)
+ \mathop{\sum_{\text{$v$ finite}}}_{v(\ideal)>0} d_v \lambda_v(P_1 - P_2)
.\]
By Lemmas \ref{lem:potgood}, \ref{lem:potmult} and \ref{lem:tetris},
together with (\ref{eq:frenchfries}),
\[\begin{aligned}
\sum_{v \in T} d_v \lambda_v(P_1 - P_2) &\geq (1-\epsilon) \sum_{v \in T} d_v 
  \min_{j=1,2}(\lambda_v(P_j))\\
&\geq (1 - \epsilon) \sum_{v\in T} d_v \max_{j=1,2}(\lambda_v(P_j)) -
\epsilon \max_{j=1,2} \sum_{v\in T} d_v \lambda_v(P_j),\\
\sum_{v \in S - T} d_v \lambda_v(P_1 - P_2)
&\geq (1- 2 \epsilon) \sum_{v \in S-T} d_v 
 \max_{j=1,2}(\lambda_v(P_j)) . \end{aligned}\]
Note that $\max_{j=1,2} \sum_{v \in T} d_v \lambda_v(P_j)
\leq \sum_{v \in T} d_v \max(\lambda_v(P_1), \lambda_v(P_2))$.
Thus
\[\sum_{v\in S} d_v \lambda_v(P_1 - P_2) \geq
(1 - 2\epsilon) \max_{j=1,2} \sum_{v\in S} d_v \lambda_v(P_j).\]
%(Write $x_v = d_v \lambda_v(P_1), y_v = d_v \lambda_v(P_2),
%I = \sum_{v \in S} d_v \lambda_v(P_1- P_2)$. 
%Then \begin{eqnarray*} S \geq (1-\epsilon) \sum_{v \in T}   \min(x_v, y_v)
%+ (1- 2 \epsilon) \sum_{v \in S-T} \max(x_v, y_v)
%\\ \geq 
% (1-\epsilon) \sum_{v \in T} \max(x_v, y_v)
%- \sum_{v \in T} |x_v - y_v| + (1-\epsilon) \sum_{v \in S-T}
% \max(x_v, y_v) - \epsilon \sum_{v \in S-T} \max(x_v, y_v)
%\\ \geq 
% (1-\epsilon) \sum_{v \in S} \max(x_v, y_v) 
%- \epsilon \max (\sum_{v \in T} x_v, \sum_{v \in T} y_v)
%- \epsilon \sum_{v \in S-T} \max(x_v, y_v) 
%\\ \geq (1-\epsilon) \sum_{v \in S} \max(x_v, y_v)
%- \epsilon \sum_{v \in T} \max(x_v, y_v) - 
%\epsilon \sum_{v \in S-T} \max(x_v, y_v)
%\\  \geq (1-2 \epsilon) \sum_{v \in S} \max(x_v, y_v)
%\end{eqnarray*}
Since $x(P_1)$, $x(P_2)$ are $S$-integers, and $S$ contains
all infinite places and all primes dividing the discriminant of $E$, we see:
\[\lambda_v(P_j) = \frac{1}{2} \log^+(|x(P_j)|_v) = 0\]
for $v\notin S$.

It remains to consider $\lambda_v(P_1 -P_2)$ for
$v$ finite, $v(\ideal)>0$. Let $\pr_v$ be the corresponding
prime ideal of $\mathscr{O}_K$, and $n_v$ its multiplicity
in $\ideal$. The point $P_1 - P_2$ is not $O$, but it is mapped
to origin when reduced modulo $\pr_v^{n_v}$. Hence $v(x(P_1-P_2))\leq -
2n_v$.
Therefore $\lambda_v(P_1 - P_2) \geq n_v e_v^{-1} \log(p_v)$,
where $p_v$ is the
rational prime lying under $v$ and $e_v$ the ramification
degree of $K_v$ over $\Q_p$. We note that $\sum_{v(\ideal) > 0 } d_v n_v e_v^{-1}
\log(p_v) =\log(\Norm \ideal)$ to conclude.
\end{proof}

\subsection{Slicing and packing}
We will use Prop.\ \ref{prop:affreuse}
to give an upper bound on the number of $S$-integral points
on the curve $E:y^2 = x^3 + d$. Any application of quasi-orthogonality leads fairly
naturally to a bound of the form
\begin{equation}\label{eq:troglo}
\text{$\#$ of integer points on $E$} \ll C^{\rank(E) + \# S} ,\end{equation}
for some constant $C$ (vd. \cite{GS}). However, in applications such
as estimating the size of $3$-torsion, the size of $C$ is crucial; if $C$ is
too large, one does not recover even the trivial bounds on $3$-torsion. 
One may optimize the bound by applying sphere packing (cf.\ \cite{He});
in order to make this approach particularly effective, we first
partition the set of integer points on $E$,
and then apply sphere-packing bounds to each part separately. 

For $\vec{x} = (x_i)_{1 \leq i \leq n} \in \R^n$, we set $|\vec{x}|_{1} = \sum_{1 \leq i \leq n} |x_i|$. 
Let $d(\vec{x},\vec{y}) = |\vec{x} - \vec{y}|_{1}$ be the associated metric on $\R^n$. 
Let $\R_{\geq 0} = \{z \in \R: z \geq 0\}$. 
\begin{lem}\label{lem:covak}
Let $c_1,c_2$ be positive real numbers, $0 < \epsilon < 1/2$, 
$n$ a non-negative integer. 
Let $S = \{\vec{x} \in \R_{\geq 0}^n: c_1 \leq \sum_{i=1}^{n}  x_i < c_2\}$. 
Then there is
a subset $T \subset \R_{\geq 0}^{n}$ and an explicit, absolute constant $C>0$  
such that \begin{equation}\label{eq:ubound}\# T \leq
C^n \epsilon^{-(n+1)} (1 + \log(c_2/c_1)) \end{equation}
and the $\ell_1$-balls $B(P, \epsilon |P|_{1})$,
for $P \in T$, cover all of $S$.  
\end{lem}
\begin{proof}
The idea is to slice $S$ into a union of regions where $|x|_{1}$
is almost constant and then to consider the points on a lattice in each of
these regions. 

Since we may replace $\epsilon$ by $\epsilon/2$, it suffices to
cover $S$ by balls $B(P, 2\epsilon |P|_{1})$. 
Let $\lfloor z \rfloor$ be the largest integer no greater than $z$.
For $\vec{x} =(x_i)_{1 \leq i \leq n} \in \R^n$, we set $\lfloor \vec{x} 
\rfloor = (\lfloor x_i \rfloor)_{1 \leq i \leq n}$. 
Define $M =  \frac{\log(c_2/c_1)}
 {\log(1 + \epsilon)}$. 
Set \[T = \bigcup_{0\leq m<M} \frac{c_1 \epsilon 
(1+\epsilon)^{m}}{n} 
\{\vec{y} \in \Z^{n}_{\geq 0} :
n (\epsilon^{-1}-1) \leq |\vec{y}|_{1} < n (1 + \epsilon^{-1}) \} .\]
Then $T$ has the required property: given $\vec{x} \in S$, 
set
\[m(\vec{x}) = \left\lfloor \frac{\log(|\vec{x}|_{1}/c_1)}{\log(1+\epsilon)}\right\rfloor,
\;\;\;\;
\vec{y}(\vec{x}) = \left\lfloor \frac{n  \vec{x}}{c_1 \epsilon (1 +\epsilon)^{m(\vec{x})}}  \right\rfloor .\] 
Then $P =
\frac{c_1 \epsilon (1+ \epsilon)^{m(\vec{x})}}{n}
\vec{y}(\vec{x})$ belongs to $T$. Moreover, $P$
 satisfies $d(\vec{x},P) \leq \frac{\epsilon}{1-\epsilon} |P|_{1}
\leq 2 \epsilon |P|_{1} $, by virtue of 
the fact that $d(\vec{z},
\lfloor \vec{z} \rfloor) \leq n$ for any $\vec{z} \in \R^n$. 
It remains to estimate $\# T$:

\[\begin{aligned}  \# T
&\leq \left(1 + \frac{\log(c_2/c_1)}{\log(1+\epsilon)}\right) \cdot  
\#\{\vec{y} \in \Z_{\geq 0}^n : n(\epsilon^{-1}-1) \leq |\vec{y}|_{1} 
\leq n(1+\epsilon^{-1})\} \\ 
&\leq \left(1 + \frac{\log(c_2/c_1)}{\log(1+\epsilon)}\right)  
\frac{( n(1+\epsilon^{-1}) + n)^n }{n!} .
\end{aligned}\]
The result follows by Stirling's formula. 
\end{proof}

We will need lower bounds on the canonical height. Note that there are
strong bounds for the {\em number} of points
of moderately low height \cite{Da}; 
such bounds could be used in place of the following proposition.
\begin{prop}\label{prop:lowbound}
 Let $E$ be an elliptic curve over a number field $K$.
 There is an absolute constant $0 < \kappa < 1$ such that, for every non-torsion
 point $P\in E(K)$,
 \[\hat{h}(P) > \kappa^{m + \lbrack K:\Q \rbrack} \max(1,h(j)),\]
where $m$ is the number of places of $K$ where $E$ has potentially
multiplicative
reduction, and $j = j(E)$ is the $j$-invariant of $E$. 
\end{prop}
 \begin{proof}
 By the proof of the Theorem in \cite{Sil5}, \S 4; see also \cite{Sil6},
 Thm.\ 7.\end{proof}
%A result of David \cite{Da} shows, notations being as in Proposition
%\ref{prop:lowbound},
%that the number of $P$ with $\hat{h}(P) \leq \frac{\max(1,h(j))}{[K:\Q]
%+ m}$ is $\ll (m + [K:\Q]) (1 + \frac{\log([K:\Q])}{\max(1,h(j))}) $. 

We shall use the remarkable bounds of Kabatiansky
and Levenshtein.
\begin{prop}\label{prop:kl}
Let $A(n,\theta)$ be the maximal number of points that can be arranged on
the unit sphere of
$\mathbb{R}^n$ such that the angle $\angle P_1 O P_2$ between any two of them and the origin is no smaller than $\theta$. Then
for $0 < \theta < \pi/2$, 
\begin{equation} \label{eq:packing}\frac{1}{n} \log_2 A(n,\theta) \leq
\frac{1 + \sin \theta}{2 \sin \theta} \log_2 \frac{1 + \sin \theta}{2 \sin \theta} - 
\frac{1 - \sin \theta}{2 \sin \theta} \log_2 \frac{1 - \sin \theta}{2 \sin
\theta} + o(1) ,\end{equation}
where the convergence of $o(1)\to 0$ as $n\to \infty$ is uniform and explicit 
for $\theta$ within any closed subinterval of $(0,\pi/2)$.
In particular, for $\theta = \pi/3$, we have
\[\frac{1}{n} \log_2 A(n,\theta) \leq 0.40141\dotsc .\]
\end{prop}
%WAIT! Large enough depends on what?
%60 degrees gives 0.40141...
%Is that estimate with - 0.0099... at the end any better?
\begin{proof}
See \cite{KL}; vd. also the expositions in \cite{Le} and \cite{CS}, Ch.\ 9.
\end{proof}
\begin{Rem}
%If $\theta < \pi/2$ is fixed, then $A(n,\theta)$ grows
%exponentially with $n$. 
%In our applications, $n$ will be the Mordell-Weil rank of some elliptic
%curve and we will not have good upper bounds for it. 
%We therefore wish to only invoke (\ref{eq:packing}) when the right hand side is fairly small,
%i.e. when $\theta$ is near $\pi/2$.  
For fixed $\theta > \pi/2$,
the function $A(n, \theta)$ is bounded above {\em independently} of $n$:
given $k$ unit vectors $v_1, v_2,\dots, v_k$ separated by angles of $\theta$ 
or more, 
\begin{equation}\label{eqn:unif}\begin{aligned}
0 &\leq \langle v_1 + \dots + v_k, v_1 + \dots + v_k \rangle
\leq k + k(k-1) \max_{i \neq j} \langle v_i, v_j \rangle \\ 
&\leq k + k (k-1) \cos(\theta) .\end{aligned}\end{equation}

It may hence 
not be surprising that the derivative of the right side of 
(\ref{eq:packing}) is zero for $\theta = \pi/2$. This qualitative feature
is, in fact, the crucial ingredient in our bound on $3$-torsion.
In our application, we will apply (\ref{eq:packing})
with a $\theta$ that we will have some freedom in choosing.
As $\theta$ decreases, the increase in the right-hand side of (\ref{eq:packing})
will be offset by a decrease in ``cost'' linear in $\theta$.
In the neighborhood of $\pi/2$, therefore, it will always
be advantageous to decrease $\theta$ slightly. 

We shall put this idea in practice in the following way. 
We shall partition the set of integral points on an elliptic curve
so that any two points $P,Q$ in the same part
are separated by an angle of at least $\theta$ in the Mordell-Weil 
lattice. We will then apply (\ref{eq:packing}) to bound the number of points in each
part. (We can do the same for rational points on curves
of higher genus; see Section \ref{sec:beyond}.) 
The bounds that correspond to $\theta = \pi/2$ will correspond
(at least in cases where one can bound the difference between
canonical and naive heights) to the ``uniform'' bounds of Bombieri-Pila
and Heath-Brown. Reducing $\theta$ slightly, under favorable circumstances,
gives an improvement.

The agreement between the output of this method and
the results of \cite{BP} and \cite{HBR} is no coincidence: see
the remarks after Theorem \ref{thm:astarte}. 
\end{Rem}
\subsection{Bounding integral points}\label{subs:bipbip}
In the theorem that follows, the reader might wish
to ignore the dependence on $S$ in a first reading.
The theorem asserts, in approximate language, that the number of points in $E(K,S)$ of height up to $h_0$
is bounded above by $e^{t [K:\Q] h_0 + (\beta(t)+\epsilon)r}$, where $r$ is the Mordell-Weil rank.
Here $t \in [0,1]$ is a free parameter that will be optimized in applications: the basic idea is that if $r$ is small compared to $h_0$ it is advantageous to take $t$ small, whereas in applications where
$r$ might be very large, we take $t$ close to $1$. This optimization process is formalized in Cor.\ \ref{cor:bound}. 

Roughly speaking, the proof of the Theorem proceeds, in words, as follows. We partition the points of $E(K,S)$ into points mod $\mathscr{I}$, where $\mathscr{I}$ is a suitable ideal in $\order_K$
with norm about $e^{t [K:\Q] h_0}$. Prop.\ \ref{prop:affreuse} shows that -- after some slight refinement
of this partition -- the points belonging to the same part are very well-separated in the Mordell-Weil lattice. We then apply sphere packing bounds in the form of Prop.\ \ref{prop:kl} to each part separately. The term $e^{t [K:\Q] h_0}$ arises from the number of parts, whereas the term $e^{\beta(t) r}$ arises
from the sphere packing bounds applied to each part. We finally 
note that the purpose of
 most of the auxiliary Lemmas on previous pages is to help us carry out 
the ``slight refinement'' mentioned above.  

\begin{thm}\label{thm:astarte}
Let $E$ be an elliptic curve over a number field $K$
defined by a Weierstrass equation (\ref{eq:weierstrass}). Let $S$
be a finite set of places of $K$, including all infinite places and 
all primes dividing the discriminant of $E$. 

Then, for every $h_0\geq 1$ and every choice of $t\in \lbrack 0,1\rbrack$,
 the number of $S$-integer
points of $E(K)$ of canonical height up to $h_0$ is at most  
\begin{equation} \label{eq:canbound} O_{\epsilon, \lbrack K : \mathbb{Q}\rbrack} \left( 
C^s \epsilon^{-2 (s + \lbrack K:\mathbb{Q}\rbrack)} s^{\lbrack K :
\mathbb{Q}\rbrack}  (1+\log h_0)^2
e^{t [K:\Q]  h_0 +  (\beta(t) + \epsilon) r}\right),\end{equation}
for every sufficiently small $\epsilon$,
where $r$ is the rank of $E(K)$ as a $\mathbb{Z}$-lattice, 
$s$ is $\# S$, $C$ is an absolute constant,
\begin{equation}\label{eq:dotty}\begin{aligned}
\beta(t) &= \frac{1 + f(t)}{2 f(t)} 
\log \frac{1 + f(t)}{2 f(t)}  - \frac{1 - f(t)}{2 f(t)} 
\log \frac{1 - f(t)}{2 f(t)},\\
f(t) &= \frac{\sqrt{(1+t) (3-t)}}{2} .\end{aligned}\end{equation}
for $t\in \lbrack 0,1)$. We set $\beta(1)=0$.
\end{thm}
\begin{proof} 
We first carry out a very mild partitioning (i.e., into very few parts) of $E(K,S) $
so as any two points in the same part have comparable canonical height. 
Applying Prop.\ \ref{prop:lowbound}, we see that one can cover the set $\{P \in E(K,S):
\hat{h}(P) \leq h_0\}$ by
by $\ll \epsilon^{-1} (\log(h_0) + s) $ sets of the form $\{P \in E(K,S): h_i \geq \hat{h}(P) \geq (1- \epsilon) h_i\}$.
It therefore suffices to prove the bound \eqref{eq:canbound}, with $(1+\log h_0)^2$ replaced
by $(1+\log h_0)$, just for the set of points $P$ 
satisfying $(1-\epsilon) h_0 \leq \hat{h}(P) \leq h_0$.

Suppose first that $t \neq 0$. 
Let $\overline{S}$ be the set of places of $\Q$ below $S$.
If $X = \max(\lceil e^{t h_0}\rceil, 
 (\# \overline{S})^{1 + 1/\lbrack K : \mathbb{Q}\rbrack},C_{[K:\Q]})$,
where $C_{\lbrack K : \mathbb{Q}\rbrack}$ 
is an appropiately chosen constant, there is a prime
$p$ of $\mathbb{Q}$ with $X\leq p\leq 2 X$ and $p\notin \overline{S}$.
The ideal $\mathscr{I}$ of $\mathscr{O}_K$ generated by $p$ satisfies
\begin{equation} \label{eq:num1}
\frac{\log N \mathscr{I}}{\lbrack K:\mathbb{Q}\rbrack} \geq t h_0,\;\;\;
N \mathscr{I} \ll_{[K:\Q]}
%(2 C_{\lbrack K : \mathbb{Q}} 
% \lbrack K : \mathbb{Q}\rbrack)^{2 \lbrack K: \mathbb{Q}\rbrack} 
s^{\lbrack K : \mathbb{Q}\rbrack + 1}
e^{t \lbrack K:\mathbb{Q}\rbrack h_0}
.\end{equation}
The $S$-integer points of $E(K)$ fall into at most
$O_{[K:\Q]}(\Norm \ideal)$ classes under reduction modulo $\ideal$.

Let $R$ be the set of all infinite places and all places of potentially
 multiplicative
reduction. 
For every $v\in R$, partition $E(K_v)$ into $n_v+1$ subsets, 
where
 $n_v$ is as in (\ref{eq:marmite}) for $v$ finite, and $n_v$ is as in (\ref{eq:ago}) for $v$ infinite, in both cases with $\epsilon/2$
 instead of $\epsilon$. Consider any
 tuples $(a_v)_{v\in R}$, $(b_v)_{v\in R}$ with
 $0\leq a_v\leq n_v$, $b_v\in \{0,1\}$.
 Define
 $\mathscr{B}$ to be the set of non-torsion points $P\in E(K,S)$ such that, for
 all $v\in R$, (a) $P\in W_{v,a_v}$, (b)
 $\lambda_v(P)\geq 0$ if and only if $b_v=1$. 
 We will show how to bound the cardinality
 of 
 $\mathscr{B}_{h_0} = \{P\in \mathscr{B}: (1-\epsilon) h_0  \leq \hat{h}(P)\leq h_0\}$. Combining
this with the fact that the number of sets $\mathscr{B}$ is at most
 \begin{equation} \label{eq:num2}c_0^{s} |\log \epsilon|^{s + \lbrack K:\mathbb{Q}\rbrack}
 \epsilon^{-2 \lbrack K:\mathbb{Q}\rbrack}\end{equation}
will yield the conclusion.

Let $M = (S-R) \cup \{v\in R: b_v=1\}$.
Let $l:\mathscr{B}\to \mathbb{R}_{\geq 0}^M$ be the map defined by
\[P\mapsto (d_v \lambda_v(P))_{v\in M} .\]
Since $\lambda_v(P)<0$ for $v\in S-M$, Prop.\ \ref{prop:lowbound}
implies \[|l(P)|_{1}> [K:\Q] \kappa^{s} \max(1,h(j)).\] On the other hand,
by \cite{GS}, Prop.\ 3, (1), we have that
$\sum_{v \notin M}d_v \lambda_v(P) \geq -\frac{1}{24} h_K(j) - 3 \lbrack
K:\mathbb{Q}\rbrack$, and thus $|l(P)|_{1}\leq \lbrack
K:\mathbb{Q}\rbrack
(h_0 + 3 + h(j)/24)$ whenever
$P \in \mathscr{B}_{h_0}$. 
By Lemma \ref{lem:covak}, we can cover
$l(\mathscr{B}_{h_0})$ by at most
\begin{equation} \label{eq:num3}O(c_1^s \epsilon^{-(s+1)} \log(h_0+1)
 )\end{equation}
balls
$B(\mathbf{x},\frac{\epsilon}{8} | \mathbf{x} |_{1})$ in the metric
$|\cdot |_{1}$.
For $P_1, P_2 \in \mathscr{B}_{h_0}$ with $l(P_1), l(P_2) \in
B(\mathbf{x},\frac{\epsilon}{8} |\mathbf{x}|_{1})$, we have
$|l(P_1) - l(P_2)|_{1} \leq \frac{\epsilon}{4} |\mathbf{x}|_{1}
\leq \frac{\epsilon}{2} \max_{j=1,2} |l(P_j)|_{1}$. Suppose $P_1$
and $P_2$ have the same reduction modulo $\mathscr{I}$. Then, by Prop.\
\ref{prop:affreuse},
\begin{equation}\label{eq:argante}
\hat{h}(P_1-P_2)\geq (1-\epsilon) \max_{j=1,2} \hat{h}(P_j)
+ \frac{\log(\Norm \ideal)}{\lbrack K : \mathbb{Q}\rbrack}
\geq (1 + t - \epsilon) \max_{j=1,2} \hat{h}(P_j). \end{equation}

Embed the Mordell-Weil lattice $E(K)$ modulo torsion into $\mathbb{R}^{\rank(E)}$ so as to send $\hat{h}$ to the square of the Euclidean height. 
Since $\hat{h}(P_1), \hat{h}(P_2),
\hat{h}(P_1-P_2)>0$, the images $Q_1, Q_2\in \mathbb{R}^{\rank(E)}$ of
$P_1$ and $P_2$ are different from each other and from the origin $O$. 
By (\ref{eq:argante}), and the fact that $\hat{h}(P_1), \hat{h}(P_2)$ lie in the interval $[(1-\epsilon) h_0, h_0]$,  the angle
$\angle Q_1OQ_2$ is at least $\arccos \frac{1 - t + O(\epsilon)}{2}$.
We may now apply the KL bound (Prop.\ \ref{prop:kl}), and
obtain that there are at most 
$e^{(\beta(t)+O(\epsilon)) r}\cdot O_{\lbrack K:\mathbb{Q}\rbrack}(1)$ 
points of $\mathscr{B}_{h_0}$ with image in a given
ball $B(\mathbf{x}, \frac{\epsilon}{8} |\mathbf{x}|_{1})$
and with prescribed reduction modulo $\ideal$. (The factor
$O_{\lbrack K:\mathbb{Q}\rbrack}(1)$ is an upper bound (\cite{Me})
on the number of torsion points in $E(K)$.)
 Combining this with our
estimates for the number of possibilities for reduction
mod $\mathscr{I}$ (\ref{eq:num1}), the 
number of sets $\mathscr{B}$ (\ref{eq:num2}), 
and the number of balls $B(\mathbf{x}, \dots)$ (\ref{eq:num3}),
we obtain the statement of the Theorem. 

In the case of $t=0$, we proceed as above but without using $\ideal$. 
\end{proof}

\begin{Rem}
Note that $t=0$ gives a pure application of sphere-packing,
whereas $t=1$ recovers a bound of the quality of $c^{h_0}$ 
(for some constant $c$) with
almost no dependence on the rank.
For our bound
on $3$-torsion (Theorem \ref{thm:classbound}) we will apply the result
with $t\in (0,1)$ optimized; for the result on elliptic curves
(Theorem \ref{thm:ecurve}) we will apply it with $t=0$. 

The bound with $t=1$ is very closely related to the Bombieri-Pila
bound \cite{BP}.   
To see this, take for a moment $K = \Q$
and let $E$ be given by a Weierstrass equation (\ref{eq:weierstrass}). 
The canonical height of the integral point
$P=(x,y)$ on $E$ is given by
 $\hat{h}(P) = \frac{\log(x)}{2} + O_E(1)$; we shall ignore the term
$O_E(1)$ for the sake of exposition. If $N$ is large, then any integral point
$P= (x,y)$
on $E$ with $|x| \leq N, |y| \leq N$ has in fact $|x| \ll N^{2/3}$
and thus
$\hat{h}(P) \ll \log(N)/3$. Then the bound given by Theorem \ref{thm:astarte}
shows that the number of such points is at most $O(N^{1/3+ \epsilon})$,
which agrees with the bound of \cite{BP} in the case of degree $3$.

This apparent coincidence is a sign of a deeper
parallelism between the two methods. 
Suppose one attempts
to carry through the proof of Theorem \ref{thm:astarte} with $t > 1$. In
other words, we choose the auxiliary ideal $\ideal$ to satisfy
$\log(\Norm\ideal)
= 1.000001 [K:\Q] h_0$. In this case, the remark after Prop.\ \ref{prop:kl}
 shows that the number of integral
points on $E$ with height $\leq h_0$
and reducing to a fixed point modulo $\ideal$
is bounded independently of the rank of $E(K)$.
This is precisely what \cite{BP} and \cite{HBR} prove, as follows:
first, let $L$ be a large integer. One constructs 
a certain meromorphic function $f$ on $E^L$
such that $f$ vanishes to high order along the diagonally embedded $E$.
If $P_1, \dots, P_L$ all reduce to the same point (modulo $\ideal$)
then $(P_1, \dots, P_L) \in E^L$ is $\ideal$-adically near
the diagonal, so $f(P_1, \dots, P_L)$ must be divisible
by a high power of $\ideal$. On the other hand, its archimedean
norm is not too large; if $L$ and $\ideal$ are chosen correctly, one obtains
thus a contradiction.

Remarkably, the same function $f$ also lurks among our methods.
If one were to carry
out the proof of Theorem \ref{thm:astarte} with $t>1$ as suggested,
using (\ref{eqn:unif}) instead of sphere-packing, the crucial ingredient is
the fact that $\langle P_1 + \dots +P_L, P_1 + \dots +P_L \rangle \geq 0$
for any points $P_1, \dots, P_L \in E(K)$. This 
may be equivalently phrased:
$L \sum_{i} \langle P_i, P_i \rangle -  \sum_{\{i,j\}} \langle P_i- P_j,
P_i - P_j \rangle \geq 0$, where the latter sum is taken over unordered subsets 
$\{i,j\}$ of size $2$. 

Now the expression 
$(P_1, \dots, P_L) \mapsto  L \sum_{i=1}^{L} \langle P_i, P_i
\rangle
- \sum_{\{i,j\}} \langle P_i - P_j , P_i- P_j \rangle  $
is a Weil height on $E^L$ with respect to a certain divisor $D$.
(Denoting by $\pi_i: E^L \rightarrow E, \pi_{ij}: E^L \rightarrow E^2$
the projections onto the $i$th and $ij$th factors, for $i \neq j$, and
by $(O)$ and $\Delta$ the divisors on $E$ and $E^2$
defined by the origin and diagonal respectively, 
we can take $D = L\sum_{i} \pi_i^{*}((O))  - \sum_{\{i,j\}}
\pi_{ij}^{*} \Delta$.) 
From this point of view, the assertion that this height is always {\em
positive} is (more or less) the assertion that $D$ is {\em effective},
i.e. that there is a meromorphic function $f$ on $E^L$
such that $D+(f)  \geq 0$. 
It can be verified that, with appropriate choices, 
this function can be taken to be the function $f$ discussed above.

One can push this further to an almost word-for-word
translation from one method to another.
%for instance, the
%contradiction obtained above using $f$ translates to the positivity
%of the Weil height w.r.t. $D$. 
On the other hand,
when $t < 1$, the proof of Theorem \ref{thm:astarte} begins
to use, in an essential way, the geometry of elliptic curves -- one may say:
the geometry of curves of non-zero genus -- and the translation fails. 
This is hardly surprising, as the Bombieri--Pila bounds are often tight
for rational curves. 
\end{Rem}

\begin{defn}
We define
\begin{equation} \label{eqn:alphadef}
\alpha(x) =
\min(x t + \beta(t): 0 \leq t \leq 1)
\end{equation}
for $x\geq 0$, where $\beta$ is as in (\ref{eq:dotty}).
We set  $\alpha(\infty) = \beta(0)$.
\end{defn}
\begin{cor}\label{cor:bound}
Let $E$ be an elliptic curve over a number field $K$. Let $S$
be a finite set of places of $K$, including all infinite places and
all primes dividing the discriminant of $E$.
Let $R\geq \max(1, \rank_{\mathbb{Z}} E(\mathbb{Q}))$.
Then, for every $h_0\geq 1$,
the number of $S$-integer points of $E(K)$ of canonical height
up to $h_0$ is at most  \begin{equation} \label{eq:b} O_{\epsilon, \lbrack K : \mathbb{Q}\rbrack} \left(
C^s \epsilon^{-2 (s + \lbrack K:\mathbb{Q}\rbrack)} s^{\lbrack K :
\mathbb{Q}\rbrack}  (1+\log h_0)^2
e^{R\cdot \alpha(\frac{h_0 \lbrack K : \mathbb{Q}\rbrack}{R}) + \epsilon R }
\right)\end{equation} 
for every sufficiently small $\epsilon$, where
$s$ is $\# S$ and $C$ is an absolute constant.
\end{cor}
\begin{proof}
The statement is simply that of Thm.\ \ref{thm:astarte} with $t$
optimized.
\end{proof}

\begin{Rem}
Let $S$ to be the set of all infinite places and all primes of bad
reduction, and assume, for simplicity, that $K = \mathbb{Q}$.
Assume that $h_0> c \max(\log \Delta,h(j))$ for some constant $c$.
Then the main contribution to (\ref{eq:b}) is given by
\[
e^{R\cdot \alpha(h_0/R)}.\]
Since $\beta'(1) = 0$, the minimum of $x t + \beta(t)$ is attained
to the left of $t=1$. Since $h_0>c \log \Delta \gg R$, we actually
have $\alpha(h_0/R) < (1 - \delta_0) h_0/R$ for some constant
$\delta_0>0$ depending only on $c$. We obtain a bound of the type
\begin{equation}\label{eq:hj}
\# E(K,S) \ll e^{(1 - \delta_1) h_0}\end{equation}
for any $\delta_1<\delta_0$. As remarked after Thm.\ \ref{thm:astarte},
the bound $e^{h_0}$ would be obtained if we proceeded as in \cite{BP} and
\cite{HBR}; thus (\ref{eq:hj}) gives an improvement in the exponent.
\end{Rem}

\subsection{Quantitative consequences of bounds on the height}
There is a long tradition -- starting with \cite{Ba} -- of effective
upper bounds on the height of integral points on an elliptic curve.
It is clear that any such bound yields a quantitative result, i.e., an
effective upper bound on the number of integral points.

We will see how
upper bounds on heights can be combined with pure quasi-orthogonality
so as to show that $\# E(K,S)$ is essentially bounded by a power of the
discriminant $\Delta$ of $E$.  
There are already bounds of a comparable quality in the literature; in
particular, \cite{ES} can be used to bound $\# E(\mathbb{Q},\{\infty\})$ 
by a power of $\Delta$. What we have here is simply an improvement
in the exponent. In the next section, we will be in a situation in which
exponents are crucial; we will also be able to take advantage of complex
multiplication to reduce our exponents further. 

We note that for our purposes it is very important that the available
bounds for integral points have the property that they bound the canonical
height by a power (or at least a sub-exponential function) of the coefficients of the elliptic curve. 
In our context we will use a very strong bound due to Hajdu and Herendi \cite{HjHr}. 

In what follows we take $K = \Q$
for simplicity.

\begin{prop}\label{prop:aghyut}
Let $E$ be an elliptic curve over $\mathbb{Q}$
defined by a Weierstrass equation of the form
$y^2 = x^3 + a x + b$, where $a,b\in \mathbb{Z}$.
Let $S$ be a finite
set of places of $\mathbb{Q}$, including the infinite place. Then,
for any $P\in E(\mathbb{Q},S)$,
\[\hat{h}(P) \leq c_1 p^{c_2} (1 + \log p)^{c_3 (s+1)} (s+1)^{c_4 (s+1)} 
                             |\Delta|^{c_5} \log H_E,\]
where $s= \#S$, $p$ is the largest prime in $S$ (set $p=1$ if
$S=\infty$), $\Delta$ is the discriminant of $E$,
$H_E = \max(|a|,|b|)$, 
and  $c_1, c_2, \dotsc,c_5$ are explicit absolute constants.
\end{prop}
\begin{proof}
See \cite{HjHr}, Thm.\ 2.  See 
\cite{Si3} for bounds on $|\hat{h}(P) - \frac{1}{2} h(x(P))|$.
\end{proof}

\begin{cor}\label{cor:agobio}
Let $E$ be an elliptic curve over $\Q$ defined by
a Weierstrass equation with integer coefficients. 
Let $S$ be a finite set of places
of $\Q$, including $\infty$ and
all primes dividing the discriminant of $E$. Then the number
of $S$-integer points on $E(\Q)$ is at most
\[O_{\epsilon} \left( 
C^s \epsilon^{-2 (s + 1)} 
(\log |\Delta| + \log p)^2
e^{(\beta(0) + \epsilon) r}\right)\]
for every sufficiently small $\epsilon$,
where $r$ is the rank of $E(\Q)$ as a $\mathbb{Z}$-lattice, 
$s$ is $\# S$, $C$ is an absolute constant, $p$ is the largest prime
in $S$, $\Delta$ is the discriminant of $E$,
 and $\beta(t)$ is as in (\ref{eq:dotty}). Numerically, 
$\beta(0) = 0.2782\dotsc$.
\end{cor}
\begin{proof}
We may, without loss of generality, assume that $E$ is given
by $y^2 = x^3 + ax + b$ with $a,b \in \mathbb{Z}$. 
(Indeed, we may make a linear substitution of variables transforming $E$
to this form, carrying integer points to integer
points, and increasing $\log|\Delta|$ by at most an absolutely bounded amount.)
The result is then immediate from Thm.\ \ref{thm:astarte} and Prop.\ \ref{prop:aghyut}.
(Note that $\log \log H_E \ll \log |\Delta|$, by Prop.\ \ref{prop:aghyut}
itself applied to $y^2 = x^3 - 27\Delta$.)
\end{proof}
\begin{cor}\label{cor:newcoeur}
Let $E$ be an elliptic curve over $\Q$ defined by a Weierstrass equation
with integer coefficients. 
Then the number of integer points
on $E(\mathbb{Q})$ is at most
\[O_{\epsilon} \left( 
|\Delta|^{\frac{\beta(0)}{2 \log 2} + \epsilon}\right),\]
for every sufficiently small $\epsilon$, where $\Delta$ is the discriminant of $E$ and
$\beta(t)$ is as in (\ref{eq:dotty}). Numerically,
$\frac{\beta(0)}{2 \log 2} = 0.20070\dotsc$.
\end{cor}
\begin{proof}
We can take $E$ to be given by an equation of the form $y^2 = x^3 + a x + b$,
$a,b\in \mathbb{Z}$. Let $\epsilon_0$ be sufficiently small. 
By Cor.\ \ref{cor:agobio}, we obtain a bound
of
\[O_{\epsilon_0}\left(|\Delta|^{\epsilon_0} 
\epsilon_0^{-2 (\omega(\Delta) + 1)} \log(|\Delta|)^2
e^{(\beta(0) + \epsilon_0) r}\right),\]
where $r$ is the rank of $E(K)$. 
%Again by Prop.\ \ref{prop:pinter},
% applied this time to $y^2 = x^3 - 27\Delta$, we obtain
%$\log \log H_E \ll |\Delta|^{\epsilon_0}$.
Let $K$ be the cubic field 
generated by a root of $x^3 + a x + b = 0$. Then
$r \leq \log_2 h_2(K) + o(\log |\Delta|)$ by \cite{BK}, Prop.\ 7.1. (If
$x^3 + a x + b$ is not irreducible, a stronger bound follows by 
\cite{Ma}, Prop 9.8(b).) 
Since the discriminant of $K$ divides $\Delta$, we see that
$h_2(K)\leq h(K)\ll \Delta^{1/2+\epsilon_0}$. Finally,
$2 (\omega(\Delta)+1) < \epsilon_0 \log |\Delta|$ for
$|\Delta|$ large enough, and thus
$\epsilon_0^{-2 (\omega(\Delta)+1)}
< |\Delta|^{|\log \epsilon_0| \epsilon_0}$. We set $\epsilon_0$ small enough
in terms of $\epsilon$, and are done.
\end{proof}
%We will insert a comment about Evertse-Silverman here, as soon as we
%figure out what their bound really is.
\begin{Rem}
Corollary \ref{cor:newcoeur} improves on the bound
 $O_{\epsilon}(|\Delta|^{1/2 + \epsilon})$ proven by W. Schmidt
 (\cite{Schm}, Thm.\ 1) on the basis of the results in \cite{ES}. The
 exponent
 $1/2$ arises from the trivial bound $h_2(L)\leq h(L) \ll_{\epsilon}
 \Delta^{1/2 + \epsilon}$,
 where $L$ is a cubic field over $\mathbb{Q}$ of discriminant $\Delta$.   
 
 One of our main tasks in the following section will be to do better than
 Cor.\ \ref{cor:newcoeur} in the case of Mordell equations.
 We have not been able to improve on $h_2(L)\ll \Delta^{1/2 + \epsilon}$,
 but, as we will see, we can improve on $h_3(\mathbb{Q}(\sqrt{D})) \ll
 D^{1/2 + \epsilon}$.
  Note
 that Cor.\ \ref{cor:newcoeur} would already be enough to break current
 bounds on the number of elliptic curves of given conductor (cf.\ Thm.\
 \ref{thm:ecurve}).
\end{Rem}
\section{Elliptic curves and $3$-torsion}\label{sec:ghatyar}
Throughout this section, let $D$ be a nonzero integer. 
We denote by $E_D$ the elliptic curve $y^2 = x^3 + D$. 
Suppose, for simplicity, that $D$ is negative; as we will see, we can assume as much by Scholz's
reflection principle. 
We may bound the class number $h_3(\mathbb{Q}(\sqrt{D}))$ from above by the 
number of integer points on $E_{D \delta^2}$, $1\leq \delta \ll |D|^{1/4}$.
We then apply Cor.\ \ref{cor:bound} to bound the number of integer points in 
terms of the rank of $E_{D \delta^2}$. Since $E_{D \delta^2}$ has complex
multiplication, one may do a CM-descent and thereby bound the rank of 
$E_{D \delta^2}$ in terms of $h_3(\mathbb{Q}(\sqrt{D}))$. 

We thus establish a feedback that, once started, lowers 
$h_3(\mathbb{Q}(\sqrt{D}))$ to an equilibrium point. Note that
Thm.\ \ref{thm:astarte} with $t=0$ (or $t=1$) would be insufficient to start 
the loop; only a mixed bound will do where a pure bound will not.

The problem of counting elliptic curves of given conductor also reduces to 
counting points on curves of the form $E_{D \delta^2}$. Again, their rank
may be bounded by means of a CM-descent, and the new bounds on 
$h_3(\mathbb{Q}(\sqrt{D}))$ can thus be applied.

\subsection{Bounds for $3$-torsion and cubic fields of fixed discriminant}
\begin{lem} \label{lem:rankbound}
The rank of $E_D(\Q)$ satisfies
$$\mathrm{rank}_{\Z} E_D(\Q) \leq A + B \omega(D) + 2 \log_3 h_3(\Q(
\sqrt{D}))$$
for some absolute constants $A,B$. 
\end{lem}
\begin{proof} 
See, e.g., \cite{fouvry}, Prop.\ 2. 
\end{proof}

\begin{thm} \label{thm:classbound} Let $D$ be a positive integer. Suppose
$-D$ is a fundamental discriminant. Then, for every $\epsilon>0$,
\begin{equation} \label{eq:classbound} 
h_3(\Q(\sqrt{-D})) \ll_{\epsilon} D^{\lambda+\epsilon},\;\;\;\;h_3(\Q(\sqrt{3 D}))\ll_{\epsilon} D^{\lambda+\epsilon},\end{equation}
where $\lambda$ is the unique real solution in the range 
$\lambda \in (0.4,0.5)$ to the equation
\begin{equation} \label{eqn:lambdadef}\lambda = 1/4 + \frac{2
\lambda}{\log(3)} \alpha \left(\frac{\log(3)}{8 \lambda}\right) .\end{equation}
Numerically, $\lambda = 0.44178...$
\end{thm}
\begin{proof}
Let $\iota$ be an embedding of $\Q(\sqrt{-D})$ into $\C$. Let
$\mathfrak{a}$ be an ideal of $\mathscr{O}_{\Q(\sqrt{-D})}$.
Then $\iota(\mathfrak{a})$ is a lattice in $\C$ of covolume
$N(\mathfrak{a}) D^{1/2}$. Minkowski's theorem
(see \cite{siegel}) shows that $\iota(\mathfrak{a})$
contains $x \in \C$ with $|x| \ll N(\mathfrak{a})^{1/2} D^{1/4}$,
the implicit constant being absolute. 
This implies that $\mathfrak{a}$
contains an element $\alpha$ of norm $\ll N(\mathfrak{a}) \sqrt{D}$.
Then $\mathfrak{a}^{-1} \cdot \alpha$ 
is an integral ideal in the same ideal class as $\mathfrak{a}^{-1}$ 
but of norm $\ll D^{1/2}$. 

In particular, any $3$-torsion class in the ideal class group
of $\Q(\sqrt{-D})$ 
has a representative $\mathfrak{a}$ that satisfies $N(\mathfrak{a}) \ll D^{1/2}$. 
Since $\mathfrak{a}^3$ is principal, it follows that
$\mathfrak{a}^3 = (y + \delta \sqrt{-D})$
where $y + \delta \sqrt{-D} \in \mathscr{O}_{\Q(\sqrt{-D})}$.
Thus
$N(\mathfrak{a})^3 = y^2 + D \delta^2$. 
Since $y + \delta \sqrt{-D}$ is an integer in $\Q(\sqrt{-D})$,  
we know that $2 y$ and $2 \delta$ are integers. 

Since $N(\mathfrak{a}) \ll D^{1/2}$, the $3$-torsion class represented
by $\mathfrak{a}$ has given us a solution to
\begin{equation} \label{eq:fred}
(4x)^3 = (8y)^2  + D (8\delta)^2 \ \  (4 x, 8y, 8\delta \in \Z, |x| \ll D^{1/2}, |y|
\ll D^{3/4}, |\delta| \ll D^{1/4}).\end{equation}
Conversely, any solution to (\ref{eq:fred}) determines $\mathfrak{a}$ up to
$\ll D^{\epsilon}$ possibilities. 
We have therefore deduced that $h_3(\Q(\sqrt{-D}))$ is at most a constant
times
\begin{equation} \label{eq:fredsback}
D^{1/4 + \epsilon} \max_{|\delta| \ll D^{1/4}}
\#\{(x,y) \in E_{-D \delta^2}(\mathbb{Q},\{\infty\}): |x| \ll D^{1/2}, |y| \ll D^{3/4}\}\end{equation}

The curve $E_{- D \delta^2}$ is a twist of $E_{-1}$,
and the map $(x,y) \rightarrow (\frac{x}{D^{1/3} \delta^{2/3}}, \frac{y}{D^{1/2} \delta})$ gives an isomorphism $E_{-D \delta^2} \rightarrow E_{-1}$ over $\overline{\Q}$. 
Thus the difference
$|\hat{h}(P) - \frac{1}{2} h(\frac{x(P)}{D^{1/3} \delta^{2/3}})|$ is bounded
above by a constant
for all $P \in E(\overline{\Q})$; on the other hand
\[\frac{1}{2}  h(x(P) D^{-1/3} \delta^{-2/3})  =
\frac{1}{6} h(\frac{x(P)^3}{D\delta^{2}}) \leq
\max(\frac{1}{2} \log|x(P)|, \frac{1}{6} \log|D \delta^2|).\]
Thus any point $P=(x,y) \in E_{-D \delta^2}$ satisfying (\ref{eq:fredsback})
has $\hat{h}(P) \leq \frac{\log(D)}{4} + O(1)$. 

Let $\gamma = \limsup_{D \rightarrow \infty}
\frac{\log(h_3(\Q(\sqrt{-D})))}{\log(D)}$.
Lemma \ref{lem:rankbound} shows that, for $D$ large enough and any
$\delta \ll D^{1/4}$,
\begin{equation}\label{eq:malin}\mathrm{rank}_{\Z} E_{-D \delta^2}(\Q) 
\leq R = 
\log(D) \left(\frac{2 \gamma}{\log(3)} + o(1)\right).\end{equation}
We apply Cor.\ \ref{cor:bound}
with $S = \{\infty\} \cup \{p : p | 6 D \delta^2\}$ and $h_0 = \log(D)/4 + O(1)$,
obtaining
\begin{equation} \label{eqn:gergelyrules}
\#\{ P \in E_{-D \delta^2}(\mathbb{Q},S) : \hat{h}(P) \leq
h_0\} \ll_{\epsilon} 
D^{\frac{2 \gamma}{\log 3}
\alpha(\frac{\log(3)}{ 8 \gamma}) + \epsilon}\end{equation}
for every $\epsilon>0$. 

By (\ref{eq:fredsback}) and (\ref{eqn:gergelyrules}), we conclude that
\begin{equation}\label{eq:iter}
\gamma \leq \frac{1}{4} + \frac{2 \gamma}{\log 3} \alpha\left(
\frac{\log 3}{8 \gamma}\right).\end{equation}
One has the {\em a priori} bound $\gamma \leq 1/2$. 
We iterate  (\ref{eq:iter}). 
Apply Scholz's reflection principle (\cite{scholz}) to obtain
 $h_3(\Q(\sqrt{3 D}))\ll D^{\lambda+\epsilon}$ therefrom.
\end{proof}

\begin{cor} \label{cor:cor1}
The number of cubic extensions of $\Q$
of discriminant $D$ is $O(|D|^{\lambda+\epsilon})$, where $\lambda$ is as in
Thm.\  \ref{thm:classbound}.
\end{cor}
\begin{proof}
 This is an immediate consequence of Thm.\
\ref{thm:classbound}; see \cite{hasse}, Satz 7. 
\end{proof}

\begin{Rem} 
The best previously known bound was the trivial one, namely that
$h_3(\mathbb{Q}(\sqrt{-D})) \leq h(\mathbb{Q}(\sqrt{-D})) \ll_{\epsilon}
D^{1/2 + \epsilon}$. 
The conditional results known to the authors are as follows.
S. Wong has shown (\cite{wong}) that the Birch-Swinnerton-Dyer conjecture,
together with the Riemann hypothesis for the $L$-functions
of elliptic curves, implies that $h_3(\Q(\sqrt{-D}))
\ll_{\epsilon} D^{1/4 +\epsilon}$.

Let $\chi_D$
be the quadratic Dirichlet character
associated to $\Q(\sqrt{-D})$. Then the Riemann hypothesis for 
$L(s,\chi_D)$ alone implies $h_3(\Q(\sqrt{-D})) \ll_{\epsilon}
D^{1/3 + \epsilon}$. We 
sketch the proof communicated to us by Soundararajan;
see also the remark at the end of \cite{sound}. 

Assume $D \equiv 1 \; \mathrm{mod} \, 4$ for simplicity. 
Let $\sigma$ be the Galois automorphism of $K=\Q(\sqrt{-D})$ over $\Q$.
Assuming the Riemann hypothesis for $L(s, \chi_D)$
one sees that there are $\gg_{\epsilon} D^{1/6 -\epsilon}$
primes $p$ with $p < D^{1/6}$ and $\chi_D(p) = 1$;
equivalently, 
there are $\gg_{\epsilon} 
D^{1/6 - \epsilon}$ prime ideals $\mathfrak{p}$
of $\mathscr{O}_{K}$
with $N\mathfrak{p} < D^{1/6}$ and $N \mathfrak{p}$ prime.
If two such distinct
ideals $\mathfrak{p}_1, \mathfrak{p}_2$ represented the same class
in the quotient group $\Cl(\mathscr{O}_{K})/\Cl(\mathscr{O}_{K})[3]$, 
then $\mathfrak{p}_1^{\sigma} \mathfrak{p}_2$ would represent a $3$-torsion
ideal class; in particular, $N(\mathfrak{p}_1^{\sigma}
\mathfrak{p}_2)^3$ would be, as in the proof of Thm.\ \ref{thm:classbound},
an integer of the form $x^2 + D y^2$. Since $N(\mathfrak{p}_1^{\sigma}
\mathfrak{p}_2)^3 < D$, this forces $y=0$, leading to a 
contradiction.  This shows that 
$\#(\Cl(\mathscr{O}_K)/\Cl(\mathscr{O}_K)[3]) \gg_{\epsilon} D^{1/6 - \epsilon}$,
which gives $\#(\Cl(\mathscr{O}_K)[3])  \ll_{\epsilon} D^{1/3+\epsilon}$ as desired. 
\end{Rem}

\subsection{Elliptic curves of given conductor}
The following is well-known; see, e.g., \cite{BS}, 
pf. of Thm.\ 1.

\begin{lem}\label{lem:mylem}
Let $S$ be a set of finite places of $\mathbb{Q}$. There is a map from
the set of all isomorphism classes of elliptic curves over $\mathbb{Q}$
with good reduction outside $S \cup \{2,3\}$
to  the union $\bigcup_C E_C(\mathbb{Q}, S \cup \{2,3\})$
of the sets of $(S\cup \{2,3\})$-integer points on each of the curves
\begin{equation}\label{eq:myleq}E_C: y^2 = x^3 + C,\end{equation}
where $C = \prod_{p\in S\cup \{2,3\}} p^{a_p}$, $0\leq a_p\leq 5$.
Each fiber has size at most $2^{\#S + 3}$. 
\end{lem}
\begin{proof} 
%Let $c_4$, $c_6$ be parameters describing a given
%curve with good reduction outside $S\cup \{2,3\}$; we choose
%$c_4$, $c_6$ to be the minimal such pair of parameters with
%$c_4, c_6\in \mathbb{Z}$. Then
%$1728 \Delta = c_4^3 - c_6^2$ is divisible only by primes in $S \cup \{2,3\}$.
%Let $r^6$ be the largest sixth power dividing $1728 \Delta$.
%Let $C = - 1728 \Delta/r^6$. We have
%\[\left(\frac{c_4}{r^3}\right)^2 = \left(\frac{c_6}{r^2}\right)^3 + C,\]
%and $C$ is of the form
%$C = \prod_{p\in S\cup \{2,3\}} p^{a_p}$, $0\leq a_p\leq 5$.
%That the size of each fiber is at most $2^{\#S+3}$ is shown in
See \cite{BS}, p.\ 100.
\end{proof}

\begin{thm} \label{thm:ecurve}
The number of elliptic curves over $\Q$ of conductor $N$
is \[O_\epsilon(N^{\gamma+\epsilon})\] for every $\epsilon>0$, 
where $\gamma = \frac{2 \lambda \cdot \beta(0)}{ \log 3}$,
$\beta$ is as in Thm.\ \ref{thm:astarte} and
$\lambda$ is as in (\ref{eqn:lambdadef}); numerically, 
$\gamma = 0.22377...$. 
\end{thm}
\begin{proof}
Let $S = \{p : p|N\}$, $M = \prod_{p\in S} p$. 
In view of Lemma \ref{lem:mylem},
$$\#\{E/\Q: E \mbox{ has good reduction outside of }S\} \ll 6^{2 \# S}
\max_{C} \#E_C(\mathbb{Q},S),$$
where the maximum is taken over all
$C = \prod_{p\in S\cup \{2,3\}} p^{a_p}$, $0\leq a_p\leq 5$.

By Cor.\ \ref{cor:agobio}, Lem.\ \ref{lem:rankbound}, and Thm.\ \ref{thm:classbound}
we obtain:
\[\begin{aligned}
\max_C \# E_C(\mathbb{Q},S) 
&\ll_{\epsilon} 
 N^{\epsilon} \max_C e^{\rank(E_C(\mathbb{Q})) \cdot (\beta(0) + \epsilon)}
 \\ &\ll_{\epsilon} 
 N^{\epsilon} e^{2 
\log_3(h_3(\mathbb{Q}(\sqrt{-N}))) (\beta(0) + \epsilon)}
\ll_{\epsilon} N^{2 \beta(0) \lambda / \log 3 + \epsilon} .
\end{aligned} \]

\begin{Rem}
 The above argument shows that, on any Mordell curve $E:y^2 = x^3 + D$,
where $D$ is a
 rational integer, there are at most $O(D^{0.22377\dotsc})$ integer points. 
 Notice the improvement over
  Cor.\ \ref{cor:newcoeur}. We are using, of course, the fact that
Mordell curves have complex multiplication.
\end{Rem}

\end{proof}
\section{Rational points: beyond $2/d$} \label{sec:beyond}
The technique here is also applicable to counting rational
points on curves of genus $\geq 1$. This will be pursued
in more detail in a separate paper; here we
content ourselves with indicating, in an approximate fashion, 
how one can use the method of this paper
to bound the number of points on a curve of higher genus
without knowing the rank of its Jacobian. Recall that Heath-Brown \cite{HBR} has shown
that if $C$ is (for example) a plane irreducible curve of degree $d$,
then the number of points in $C(\Q)$ of naive height $\leq H_0$
is $O_{d,\epsilon}(H_0^{2/d+\epsilon})$; Elkies independently proved a related bound \cite{El}
with a view to algorithmic applications. 

The method of this paper, roughly speaking, 
recovers the exponent $2/d$ for curves of higher genus, provided that
we may completely ignore the factors of $O(1)$ that
arise when dealing with Weil height functions. When further simplifications
are valid, the procedure delivers an exponent lower than $2/d$.

Let $C$ be a proper smooth curve of genus $\geq 1$ over a number field $K$. 
To further simplify matters, let us assume that $C$ has a
$K$-rational point $a$. 
The factors of $O(1)$ that occur 
in the computations below depend both on $C$ and $a$.
Let
$h_a: C(\overline{K}) \rightarrow \R$ be a Weil
height with respect to the divisor $(a)$,
and suppose we are interested in bounds for the number
of points $P \in C(K)$ with $h_a(P) \leq h_0$. 
Note that here $h_a$ will denote a Weil height ``over $K$'',
not an absolute height normalized by a factor $\frac{1}{[K:\Q]}$
(see the difference
between (\ref{eqn:htrel}) and (\ref{eqn:htabs})). For complete conformity
with our previous notation we should denote it $h_{a,K}$,
but we shall suppress the $K$ subscript for typographical ease. 

Let $J$ be the Jacobian of $C$. 
Let $j_a : C \rightarrow J$ be the embedding
that sends $P \in C$ to $P- a \in \mathrm{Pic}^{0}(C)$,
and let $\Theta$ be the associated theta-divisor,
i.e., $j_a(C) + j_a(C) + \dots +j_a(C)$, taken $g-1$ times.
Let $\Delta \in C \times C$ be the diagonal, and $h_{\Delta}
:C(\overline{K}) \times C(\overline{K}) \rightarrow \R$
an associated Weil height.

We denote by $\langle \cdot, \cdot \rangle_{\Theta}$
the inner product on $J(\overline{K})$ induced by the canonical
height associated to $\Theta$, which agrees up to $O(1)$
with the Weil height associated to the symmetric divisor
$\frac{1}{2}(\Theta + [-1]^{*} \Theta) \in \mathrm{Pic}(C) \otimes_{\Z} \R$.
(Cf.\ \cite{hs}, B.5; the map $z \rightarrow \langle z,z \rangle_{\Theta}$
on $J(\overline{K})$ is the map $\hat{q}_{J, \Theta}$ in the notation of
\cite{hs}, Thm.\ B.5.6.)

We set $||z||_{\Theta}^2 = \langle z,z \rangle_{\Theta}$
for $z \in J(\overline{K})$. 
Note that for $x \in C(\overline{K})$
\begin{equation} \label{eqn:twoheights}||j_a(x)||_{\Theta}^2 = g h_a(x) +
O(\sqrt{1+h_a(x)})\end{equation}
(this follows from \cite{hs}, Thm.\ B.5.9,
since $j_a^{*} \Theta$ is algebraically equivalent to $g \cdot (a)$,
cf.\  \cite{hs}, Thm.\ A.8.2.1).
By the proof of Mumford's gap principle
(see \cite{hs}, Thms.\ A.8.2.1 and
B.6.5), one has, for any $x,y \in C(\overline{K})$,
\begin{equation}  \label{eqn:mum}2 \langle j_a(x), j_a(y) \rangle_{\Theta}
= h_a(x) + h_a(y)
%+ O(\sqrt{1+h_a(x)} + \sqrt{1+h_a(y)})
- h_{\Delta}(x,y) + O(1)\end{equation}

Now suppose that $x,y \in C(K)$ are chosen so that $h_a(x), h_a(y) \leq h_0$. 
The theory of local heights (\cite{serre}) shows that,
if $x \neq y \in C(K)$ reduce to the same point
modulo $\mathfrak{p}$, any prime ideal of $\mathscr{O}_K$,
then $h_{\Delta}(x,y) \geq \log(\Norm \mathfrak{p}) + O(1)$.
Indeed, the hypothesis guarantees that $(x,y)$
is $\mathfrak{p}$-adically close to the diagonal, which
forces $h_{\Delta}(x,y)$ to be large. 

For such $x,y$, (\ref{eqn:mum}) yields 
$2 \langle j_a(x), j_a(y) \rangle_{\Theta} \leq 2 h_0 -
\log(\Norm \mathfrak{p}) + O(1)$. 
Thus, if we choose $\mathfrak{p}$ so that $\log(\Norm \mathfrak{p}) > 2(1 +
\epsilon) h_0$, we have $\langle j_a(x), j_a(y) \rangle_{\Theta}
\leq -\epsilon h_0 + O(1)$. 

On the other hand, in view of (\ref{eqn:twoheights}), we have
$\max(||j_a(x)||^2_{\Theta}, ||j_a(y)||^2_{\Theta}) < gh_0 + O(1 + \sqrt{h_0})
$.
The angle $\theta_{xy}$ between the points
$j_a(x), j_a(y)$ in (the Mordell-Weil lattice of) $J(K)$ thus
satisfies
\begin{equation} \label{eq:angles} \cos(\theta_{xy}) \leq \frac{
- \epsilon h_0 +  O(1)}{g
h_0 + O(1+\sqrt{h_0})}\end{equation}

Now, for the sake of the exposition, let us ignore the factors $O(1)$
and $O(\sqrt{1+h_0})$ in (\ref{eq:angles}).
It then follows that, if $\log(\Norm \mathfrak{p}) = 2(1 + \epsilon) h_0$,
then $\cos(\theta_{xy}) \leq - \epsilon/g$. 

For reasons outlined in the remark following Proposition \ref{prop:kl}, 
the number of vectors in $\mathbb{R}^N$ all of whose
mutual angles satisfy $\cos(\theta) \leq - \epsilon$
is bounded by $O_{\epsilon}(1)$. It follows
that the number of points $P \in C(K)$ of height $h_a(P) \leq h_0$
that reduce to a fixed point in $C(\mathscr{O}_K/\mathfrak{p})$
is $O_{\epsilon}(1)$; in particular,
the number of points $P \in C(K)$ of height $h_a(P) \leq h_0$
is $\ll_{\epsilon} \Norm \mathfrak{p} = \exp(2(1 + \epsilon) h_0)$.

To recognize the exponent, note that if $C$ is a curve
of degree $d$ in a projective space $\mathbb{P}^n$,
then the naive (exponential) height $H_{\mathbb{P}^n}$ on $C$
satisfies $\log H_{\mathbb{P}^n} - d h_a = O(1 + \sqrt{h_a})$.
Thus, the number of points $P \in C(K)$
with $H_{\mathbb{P}^n} \leq H_0$ is $\ll H_0^{2/d + \epsilon}$,
recovering Heath-Brown's result. 

Further, one can 
``perturb'' this method by decreasing $\Norm \mathfrak{p}$,
as was carried out in the text for integral points on elliptic curves;
a small enough perturbation improves the exponent $2/d$.
This has been carried out in \cite{EV}, which incorporates
also some different ideas stemming from the work of Heath-Brown \cite{HBR}.

\end{document}